\documentclass[a4paper,12pt]{article}

\usepackage{amsmath}
\usepackage{amsthm}
\usepackage{amssymb}
\usepackage{amscd}
\usepackage{hyperref}

\title{Bochner-Weitzenb\"ock formulas and curvature actions on Riemannian manifolds}
\author{Yasushi Homma}
 \date{}

\theoremstyle{plain}
\newtheorem{thm}{Theorem}[section]
\newtheorem{prop}[thm]{Proposition}
\newtheorem{cor}[thm]{Corollary}
\newtheorem{lem}[thm]{Lemma}
\theoremstyle{definition}
\newtheorem{defini}[thm]{Definition}

\numberwithin{equation}{section}

\theoremstyle{remark}
\newtheorem{rem}{Remark}[section]
\newtheorem{ex}{Example}[section]

\newcommand{\la}{\langle}
\newcommand{\ra}{\rangle}

\newcommand{\Hom}{\mathrm{Hom}}
\newcommand{\End}{\mathrm{End}}
\newcommand{\id}{\mathrm{id}}

\newcommand{\pf}{\mathrm{pf}}

\allowdisplaybreaks[3]

\begin{document}
%%%%%%%%%%%%%%%%%%%%%%%%%%%%%
%%% double space %%%%%%%%%%%%
%%%%%%%%%%%%%%%%%%%%%%%%%%%%%
%\setlength{\baselineskip}{30pt}
%%%%%%%%%%%%%%%%%%%%%%%%%

\maketitle

%%%%%%%%%%%%%%%%%%%%%%%%%%%%%%
%        abstract            %
%%%%%%%%%%%%%%%%%%%%%%%%%%%%%%
\begin{abstract}
Gradients are natural first order differential operators depending on Riemannian metrics. The principal symbols of them are related to the enveloping algebra and higher Casimir elements. We give certain relations in the enveloping algebra, which induce not only identities for higher Casimir elements but also all Bochner-Weitzenb\"ock formulas for gradients. As applications, we give some vanishing theorems. 
\end{abstract}

%%%%%%%%%%%%%%%%%%%%%%%%%%%%%%%%
%%            1              %%
%%%%%%%%%%%%%%%%%%%%%%%%%%%%%%%%
\tableofcontents
\section{Introduction}\label{sec:1}
The Dirac operator is an important tool in a wide range of mathematics and physics, which is a conformally covariant first order differential operator on a spin manifold. The principal symbols are known as the Clifford multiplications and constitute the Clifford algebra. The algebra gives a lot of features of the Dirac operator and has importance itself. We can generalize the Dirac operator under the condition that operators are first order and conformally covariant \cite{F}, \cite{SW}. The operators are called \textit{gradients} or \textit{Stein-Weiss operators}, and include basic operators in Riemannian and spin geometry; the exterior derivative, the interior derivative, the conformal Killing operator, the twistor operator and the Rarita-Schwinger operator. From recent research by T. Branson et al., we know various properties and applications of gradients, which give a new direction to geometry and analysis. 

Let $(M,g)$ be an $n$-dimensional oriented Riemannian manifold and $\mathbf{SO}(M)$ be the oriented frame bundle on $M$ with structure group $SO(n)$. An irreducible unitary $SO(n)$-module $V_{\rho}$ with highest weight $\rho$ gives an associated vector bundle $\mathbf{S}_{\rho}:=\mathbf{SO}(M)\times_{\pi_{\rho}}V_{\rho}$ on $M$. The Levi-Civita connection induces a covariant derivative $\nabla$ on $\mathbf{S}_{\rho}$. We decompose the target bundle of $\nabla$ with respect to $SO(n)$ as $\mathbf{S}_{\rho}\otimes T^{\ast}_{\mathbb{C}}(M)=\oplus_{1\le i\le N}\mathbf{S}_{\lambda_i}$. Then $\nabla$ splits along the bundle decomposition, and each summand is a first order differential operator $D^{\rho}_{\lambda_i}:\Gamma(M,\mathbf{S}_{\rho})\to \Gamma(M,\mathbf{S}_{\lambda_i})$. We call this operator \textit{gradient}. In \cite{Br2}, T. Branson discusses ellipticities of gradients by using the spectral resolutions on standard spheres. As a corollary, he shows that suitable linear combinations of gradients are bundle endomorphisms depending on the Riemannian curvature, 
\begin{equation}
\sum_{1\le i\le N} b_{\lambda_i}(D^{\rho}_{\lambda_i})^{\ast}D^{\rho}_{\lambda_i}=\textrm{curvature action}.\label{eqn:1-1}
\end{equation}
These formulas are called \textit{optimal Bochner-Weitzenb\"ock formulas}. The vector $(b_{\lambda_1},\cdots, b_{\lambda_N})$ in \eqref{eqn:1-1} needs to be a solution to a certain system of linear equations (see Theorem 5.10 in \cite{Br2}). Since the solutions constitute a $[N/2]$-dimensional subspace in $\mathbb{R}^N$, there are $[N/2]$ independent identities as \eqref{eqn:1-1}. This result induces some vanishing theorems and eigenvalue estimates \cite{BH1}, \cite{BH2}, \cite{BH3}. But it would be complicated to calculate the coefficients $\{b_{\lambda_i}\}_i$ and curvature actions in \eqref{eqn:1-1}. On the other hand, D. Calderbank, P. Gauduchon and M. Herzlich discuss refined Kato inequalities and apply them to Riemannian and spin geometry \cite{CGH}. Their key observation is that the principal symbols of gradients are related to the enveloping algebra of $\mathfrak{so}(n)$, especially higher Casimir elements. In this paper, we develop the observation further, and connect higher Casimir elements to Bochner-Weitzenb\"ock formulas. Then we obtain a universal and direct construction of the coefficients $\{b_{\lambda_i}\}_i$ and curvature actions in \eqref{eqn:1-1}. Thus Bochner-Weitzenb\"ock formulas on Riemannian manifolds are induced from invariants for $\mathfrak{so}(n)$. Moreover, our construction makes it possible to give various vanishing theorems and eigenvalue estimates. 

In Section \ref{sec:2}, we give a short review to representations of $\mathfrak{so}(n)$ and set up notation. In Section \ref{sec:3}, we discuss the enveloping algebra and Casimir elements of $\mathfrak{so}(n)$. We give the universal Bochner-Weitzenb\"ock formulas (Theorem \ref{thm:3-4}). The formulas induce some identities for higher Casimir elements (Corollary \ref{cor:3-5}). In Section \ref{sec:4}, we discuss the principal symbols of gradients called \textit{Clifford homomorphisms}. We relate them with the enveloping algebra by using conformal weights. From the universal Bochner-Weitzenb\"ock formulas, we have relations for Clifford homomorphisms corresponding to \eqref{eqn:1-1} on the symbol level. We also compute eigenvalues of Casimir elements. In Section \ref{sec:5}, we define first order geometric differential operators called gradients and give some fundamental properties of them. In Section \ref{sec:6},  we define curvature endomorphisms corresponding to curvature actions in \eqref{eqn:1-1}. An interesting observation is that the curvature endomorphism associated to the Pfaffian element depends on only the conformal Weyl tensor and the scalar curvature (Proposition \ref{prop:6-3}). In Section \ref{sec:7}, we give Bochner-Weitzenb\"ock formulas or Bochner identities for gradients (Theorem \ref{thm:7-1}) and discuss linear independence of them. As a result, our formulas give the coefficients $\{b_{\lambda_i}\}_i$ and curvature actions in \eqref{eqn:1-1} explicitly. In Section \ref{sec:8}, we mention fundamental examples, gradients on spinors and differential forms. We generalize the examples to other gradients and give some vanishing theorems. In particular, by using Bochner-Weitzenb\"ock formula associated to the Pfaffian element, we give a vanishing theorem on even dimensional conformally flat manifolds. In the last section, we discuss the four dimensional case. Because of the decomposition $\Lambda^2=\Lambda^2_+\oplus \Lambda^2_-$, the curvature endomorphisms and Bochner-Weitzenb\"ock formulas split. 

 An outline of some results has been presented in a short paper \cite{H5} by the author. We discuss the details and develop them in the present paper. Then new results and many examples are presented.

%%%%%%%%%%%%%%%%%%%%%%%%%%%%%%%%
%%            2               %%
%%%%%%%%%%%%%%%%%%%%%%%%%%%%%%%%
\section{Representations of $SO(n)$ and $Spin(n)$}\label{sec:2}
We give a short review to representation theory of the special orthogonal group $SO(n)$ or the spin group $Spin(n)$ \cite{K}, \cite{Z}. Let $\mathbb{R}^n$ be the $n$-dimensional Euclidean space with inner product $\la\cdot, \cdot\ra$, and $\{e_i\}_{i=1}^n$ be its standard basis. Associating $\xi\wedge \eta$ with a skew-symmetric endomorphism
\begin{equation}
(\xi\wedge \eta)(u)=\la \xi,u\ra \eta-\la \eta,u\ra\xi \quad \textrm{for $u\in \mathbb{R}^n$}, \nonumber
\end{equation}
we identify the space of $2$-forms $\Lambda^2(\mathbb{R}^n)$ with the Lie algebra $\mathfrak{so}(n)$ of $SO(n)$ or $Spin(n)$. We set $e_{ij}:=e_i \wedge e_j$ and know that $\{e_{ij}\}_{1\le i,j\le n}$ satisfy 
\begin{equation}
\begin{split}
&e_{ij}=-e_{ji}, \quad e_{ii}=0, \\
&[e_{kl},e_{ij}]=\delta_{ki}e_{lj}+\delta_{kj}e_{il}-\delta_{il}e_{kj}-\delta_{lj}e_{ik}, \end{split}\label{eqn:2-3}
\end{equation}
and $\{e_{ij}|1\le i<j\le n\}$ constitute a basis of $\mathfrak{so}(n)$.

We choose a Cartan subalgebra $\mathfrak{h}$ of $\mathfrak{so}(n)$ with basis $\{e_{2i-1,2i}|1\le i\le m=[n/2]\}$ and fix a basis $\{\mu_i\}_{1\le i\le m}$ of $(\sqrt{-1}\mathfrak{h})^{\ast}$ such that $\mu_i(-\sqrt{-1}e_{2j-1,2j})=\delta_{ij}$. The normalized Killing form of $\mathfrak{so}(n)$ is a negative definite inner product such that $\la e_{ij},e_{kl}\ra=-\delta_{ik}\delta_{jl}$ for $i<j$ and $k<l$, and induces a positive definite inner product on $(\sqrt{-1}\mathfrak{h})^{\ast}$, $\la \mu_i,\mu_j\ra=\delta_{ij}$. Thus we think of $(\sqrt{-1}\mathfrak{h})^{\ast}$ as the $m$-dimensional Euclidean space and  denote $\mu_i$ by
\begin{equation}
\mu_i=(\underbrace{0,\cdots,0}_{i-1},1,\underbrace{0,\cdots,0}_{m-i}).
 \nonumber 
\end{equation}

We consider a finite dimensional irreducible unitary representation $(\pi,V)$ of $SO(n)$ or $Spin(n)$, and decompose the representation space $V$ into simultaneous eigenspaces with respect to $\sqrt{-1}\mathfrak{h}$. Each eigenvalue $\nu$ called \textit{weight} is an integral or half-integral linear combinations of $\{\mu_i\}_i$, namely,  $\nu=\sum \nu^i\mu_i=(\nu^1,\cdots, \nu^m)$  in $\mathbb{Z}^m\cup (\mathbb{Z}+1/2)^m$. In the lexicographical order, we have a highest weight $\rho=(\rho^1,\cdots,\rho^m)$ for $(\pi,V)$ with multiplicity one. This highest weight satisfies \textit{the dominant condition},
\begin{equation}
\begin{split}
\rho^1\ge \rho^2\ge \cdots \ge\rho^{m-1}\ge|\rho^m|, \quad &\textrm{for $n=2m$},
    \\
\rho^1\ge \rho^2\ge \cdots \ge\rho^{m-1}\ge\rho^m\ge 0, \quad &\textrm{for $n=2m+1$}.
 \end{split}\nonumber
\end{equation}
Conversely, for a dominant weight $\rho$ in $\mathbb{Z}^m\cup (\mathbb{Z}+1/2)^m$, we can construct a unique irreducible unitary representation with highest weight $\rho$ up to equivalence. Therefore we denote by $(\pi_{\rho},V_{\rho})$ an irreducible representation with highest weight $\rho$ of $SO(n)$, $Spin(n)$ or $\mathfrak{so}(n)$. Note that if $\rho$ is half-integral, then $(\pi_{\rho},V_{\rho})$ does not factor through a representation of $SO(n)$. When writing dominant weights, we denote by $k_j$ a string $j$ $k$'s for $k$ in $\mathbb{Z}\cup (\mathbb{Z}+1/2)$, and sometimes omit a terminal string of zeros. For example, the highest weight of spinor representation is 
$$
((1/2)_{m-1},\pm 1/2)=(\underbrace{1/2,1/2,\cdots,1/2}_{m-1},\pm 1/2),
$$
and the one of the $p$th exterior tensor representation on $\Lambda^p(\mathbb{R}^n\otimes \mathbb{C})$ is 
$$
(1_p)=(\underbrace{1,\cdots,1}_p, \underbrace{0,\cdots,0}_{m-p}). 
$$
%%%%%%%%%%%%%%%%%%%%%%%%%%%%%%%
%%            3               %%
%%%%%%%%%%%%%%%%%%%%%%%%%%%%%%%%
\section{Enveloping algebra and higher Casimir elements}\label{sec:3}
Let  $\mathfrak{so}(n,\mathbb{C})$ be the complexification of $\mathfrak{so}(n)$ and $U(\mathfrak{so}(n,\mathbb{C}))$ be its enveloping algebra. The enveloping algebra is the quotient algebra of the tensor algebra of $\mathfrak{so}(n,\mathbb{C})$ by two-sided ideal generated by all $(X\otimes Y-Y\otimes X-[X,Y])$ for $X,Y$ in $\mathfrak{so}(n,\mathbb{C})$. Each representation $(\pi,V)$ of $\mathfrak{so}(n)$ is lifted naturally to a representation of the enveloping algebra denoted by the same notation $(\pi,V)$.  The center $\mathfrak{Z}$ of $U(\mathfrak{so}(n,\mathbb{C}))$ is characterize as the invariant subalgebra in $U(\mathfrak{so}(n,\mathbb{C}))$ under the adjoint action of $SO(n)$. We call elements in $\mathfrak{Z}$ \textit{Casimir elements}. From Schur's lemma, we know that every Casimir element is a constant on irreducible $\mathfrak{so}(n)$-module. 

We shall construct an algebraic basis of $\mathfrak{Z}$ \cite{NR}, \cite{O}, \cite{Z}. First, we define the usual Casimir element by $c_2:=\sum_{i,j} e_{ij}e_{ji}$. It is known that the eigenvalue of $c_2$ on irreducible $\mathfrak{so}(n)$-module $V_{\rho}$ is
\begin{equation}
\pi_{\rho}(c_2)=2\la \delta+\rho,\delta+\rho\ra-2\la \delta,\delta\ra=2\la \rho,\rho\ra+4\la\rho,\delta\ra, \label{eqn:3-1}
\end{equation}
where $\delta$ is half the sum of the positive roots, 
\begin{equation}
\delta=\begin{cases}
(m-1,m-2,\cdots, 1,0)  &\textrm{for $n=2m$}, \\
(m-1/2,m-3/2,\cdots, 3/2,1/2) &\textrm{for $n=2m+1$}.
\end{cases}
\nonumber 
\end{equation}
For example, on the natural representation $(\pi_{\mu_1},V_{\mu_1})=(\pi_{\mu_1},\mathbb{C}^n)$, we have $\pi_{\mu_1}(c_2)=2(n-1)$. 

Next we construct higher Casimir elements. For each non-negative integer $q$, we define an element $e_{ij}^q$ in $U(\mathfrak{so}(n,\mathbb{C}))$ by
\begin{equation}
e_{ij}^q:=
\begin{cases}
\sum_{1\le i_1,i_2,\cdots,i_{q-1} \le n} e_{ii_1}e_{i_1i_2}\cdots e_{i_{q-1}j}  & q\ge 1, \\
   \delta_{ij} & q=0.
\end{cases}\label{eqn:3-4}
\end{equation}
This $e_{ij}^q$ behaves like $e_{ij}$ under the adjoint action of $\mathfrak{so}(n)$ as follows.
\begin{lem}\label{lem:3-1}
The elements $\{e_{ij}^q|q\in \mathbb{Z}_{\ge 0},i,j=1,\cdots,n\}$ satisfy that 
\begin{gather}
[e_{kl},e_{ij}^q]=\delta_{ki}e_{lj}^q+\delta_{kj}e_{il}^q-\delta_{il}e_{kj}^q-\delta_{lj}e_{ik}^q, \label{eqn:3-5}\\
\sum_{1\le k\le n} e_{ik}^p e_{kj}^q=e_{ij}^{p+q}. \label{eqn:3-6}
\end{gather}
\end{lem}
\begin{proof}
We calculate the adjoint action of $\mathfrak{so}(n)$ on $e_{ij}^q$. From \eqref{eqn:2-3}, we have
\begin{equation}
\begin{split}
[e_{kl},e_{ij}^q] &=\sum_{ i_1,\cdots,i_{q-1}} [e_{kl}, e_{ii_1}]e_{i_1i_2} \cdots e_{i_{q-1}j} 
 +\cdots +\sum_{i_1,\cdots,i_{q-1}} e_{ii_1}e_{i_1i_2} \cdots [e_{kl},e_{i_{q-1}j}]\\
 &=\delta_{ki}e_{lj}^q+\delta_{kj}e_{il}^q-\delta_{il}e_{kj}^q-\delta_{lj}e_{ik}^q.
\end{split}\nonumber
\end{equation}
The equation \eqref{eqn:3-6} is clear from the definition \eqref{eqn:3-4}. 
\end{proof}
The equation \eqref{eqn:3-5} yields that the trace $c_q:=\sum_i e_{ii}^q$ is an invariant element. Thus we have higher Casimir elements $\{c_q\}_{q\ge 0}$. Note that $c_0$ is equal to $n$ and $c_1$ is zero. The eigenvalues of $c_q$ on irreducible $\mathfrak{so}(n)$-modules are calculated in Section \ref{sec:4}. 

In the case of $n=2m+1$, Casimir elements $\{c_q\}_q$ generate the center $\mathfrak{Z}$ algebraically. On the other hand, in the case of $n=2m$, we need another Casimir element to generate $\mathfrak{Z}$. We define $\pf$ in $\mathfrak{Z}$ by 
\begin{equation}
\pf:=\frac{1}{(\sqrt{-1})^m 2^mm!}\sum_{\sigma \in \mathfrak{S}_{2m}} \mathrm{sign}(\sigma) e_{\sigma (1)\sigma (2)}e_{\sigma (3) \sigma (4)} \cdots e_{\sigma (2m-1)\sigma (2m)}, \label{eqn:3-7}
\end{equation}
where $\mathfrak{S}_{2m}$ is the permutation group of $\{1,\cdots,2m\}$. We call the Casimir element $\pf$ \textit{the Pfaffian element}. The following facts on these Casimir elements are known \cite{NR}, \cite{O}, \cite{Z}. 
\begin{prop}\label{prop:3-2}
\begin{enumerate}
	\item In the case of $n=2m$, $\{c_2,c_4,\cdots,c_{2m-2},\pf\}$ generate $\mathfrak{Z}$ algebraically. The eigenvalue of $\pf$ on irreducible $\mathfrak{so}(2m)$-module $V_{\rho}$ is 
\begin{equation}
\pi_{\rho}(\pf)=(\rho^1+m-1)(\rho^2+m-2)\cdots (\rho^{m-1}+1)\rho^m.\label{eqn:3-8}
\end{equation}
	\item In the case of $n=2m+1$, $\{c_2,c_4,\cdots,c_{2m}\}$ generate $\mathfrak{Z}$ algebraically.
\end{enumerate}
\end{prop}

In the above proposition, we have a problem of how $c_{2q+1}$ is realized as a polynomial in $\{c_{2p}\}_p$. To answer it, we return to $e_{ij}^q$ in $U(\mathfrak{so}(n,\mathbb{C}))$ and search how $e_{ij}^q$ is related to $e_{ji}^q$. From \eqref{eqn:3-5}, we have 
\[ 
e_{ij}^{q+1}=(1-n)e_{ij}^q+\delta_{ji}c_q-e_{ji}^q-\sum_k e_{kj}^qe_{ki}.
 \]
This equation implies that $e_{ij}^q$ is a linear combination of $\{e_{ji}^p\}_{p=0}^q$ whose coefficients are Casimir elements, $e_{ij}^q=\sum_{p=0}^q a_{q,p} e_{ji}^p$. Here $\{a_{q,p}\}_{q\ge p\ge 0}$ in $\mathfrak{Z}$ are given by a recursion formula. Since the recursion formula is a little complicated, we translate $e_{ij}^q$ to another element. We define $\hat{e}_{ij}$ by
\begin{equation}
\hat{e}_{ij}:=e_{ij}+\frac{n-1}{2}\delta_{ij}, 
\nonumber 
\end{equation}
and $\hat{e}_{ij}^q$ by 
\begin{equation}
\hat{e}_{ij}^q:=
\begin{cases}
\sum_{1\le i_1,i_2,\cdots,i_{q-1} \le n} \hat{e}_{ii_1}\hat{e}_{i_1i_2}\cdots \hat{e}_{i_{q-1}j}  & q\ge 1, \\
   \delta_{ij} & q=0.
\end{cases}\nonumber 
\end{equation}
Note that $\hat{e}_{ij}^q$ is related to $e_{ij}^q$ as
\begin{equation}
\hat{e}_{ij}^q=\sum_{p=0}^q \binom{q}{p}\left(\frac{n-1}{2}\right)^{q-p}e_{ij}^p. \nonumber 
\end{equation}
We also define the translated Casimir element $\hat{c}_q$ by $\hat{c}_q:=\sum_i\hat{e}_{ii}^q$. 
%%%%%%%%%%%%%%%%%
\begin{lem}\label{lem:3-3}
The translated elements $\{\hat{e}_{ij}^q|q\in \mathbb{Z}_{\ge 0},i,j=1,\cdots,n\}$ satisfy that
\begin{gather}
[\hat{e}_{kl},\hat{e}_{ij}^q]=\delta_{ki}\hat{e}_{lj}^q+\delta_{kj}\hat{e}_{il}^q-\delta_{il}\hat{e}_{kj}^q-\delta_{lj}\hat{e}_{ik}^q, \label{eqn:3-13}\\
\sum_k \hat{e}_{ik}^p\hat{e}_{kj}^q=\hat{e}_{ij}^{p+q}, \label{eqn:3-14}\\
\hat{e}_{ij}=-\hat{e}_{ji}+(n-1)\delta_{ij}. \label{eqn:3-15}
\end{gather}
We particularly obtain  a key relation
\begin{equation}
\hat{e}_{ij}^{q+1}=\delta_{ji}\hat{c}_q-\hat{e}_{ji}^q-\sum_k\hat{e}_{kj}^q\hat{e}_{ki}. \label{eqn:3-16}
\end{equation}
\end{lem}
%%%%%%%%%%
\begin{proof}
We can prove \eqref{eqn:3-13}-\eqref{eqn:3-15} in the same way as
 Lemma \ref{lem:3-1}. It follows from \eqref{eqn:3-13}-\eqref{eqn:3-15} that 
 \begin{equation}
\begin{split}
\hat{e}_{ij}^{q+1}=&\sum_k \{[\hat{e}_{ik},\hat{e}_{kj}^q]+\hat{e}_{kj}^q\hat{e}_{ik}\} \\
=&\sum_k \{(\delta_{ik}\hat{e}_{kj}^q+\delta_{ij}\hat{e}_{kk}^q-\delta_{kk}\hat{e}_{ij}^q-\delta_{kj}\hat{e}_{ki}^q)+\hat{e}_{kj}^q(-\hat{e}_{ki}+(n-1)\delta_{ik})\}\\
=&\hat{e}_{ij}^q+\delta_{ij}\hat{c}_q-n\hat{e}_{ij}^q-\hat{e}_{ji}^q-\sum_k \hat{e}_{kj}^q\hat{e}_{ki}+(n-1)\hat{e}_{ij}^q\\
=&\delta_{ij}\hat{c}_q-\hat{e}_{ji}^q-\sum_k\hat{e}_{kj}^q\hat{e}_{ki}.
\end{split}\nonumber
  \end{equation}
  Thus we obtain \eqref{eqn:3-16}. 
  \end{proof}
%%%%%%%%%%%  
By using the above lemma, we connect $\hat{e}_{ij}^q$ with $\hat{e}_{ji}^q$. 
%%%%%%%%%%% theorem %%%%%%%%%%%%%%%%
\begin{thm}\label{thm:3-4}
The translated element $\hat{e}_{ij}^q$ is a linear combination of $\{\hat{e}_{ji}^p\}_{p=0}^q$ whose coefficients are Casimir elements,
\begin{equation}
\hat{e}_{ij}^q=(-1)^q\hat{e}_{ji}^q-\frac{1-(-1)^q}{2}\hat{e}_{ji}^{q-1}+\sum_{p=0}^{q-1}(-1)^p\hat{c}_{q-1-p}\hat{e}_{ji}^p. \label{eqn:3-17}
\end{equation}
Thus we have
\begin{gather}
\hat{e}_{ij}^{2q}=\hat{e}_{ji}^{2q}+\sum_{p=0}^{2q-1}(-1)^p\hat{c}_{2q-1-p}\hat{e}_{ji}^p, \label{eqn:3-18} \\
\hat{e}_{ij}^{2q+1}=-\hat{e}_{ji}^{2q+1}-\hat{e}_{ji}^{2q}+\sum_{p=0}^{2q}(-1)^p\hat{c}_{2q-p}\hat{e}_{ji}^p. \label{eqn:3-19}
\end{gather}
\end{thm}
%%%%%%%%%%%
\begin{proof}
Setting $\hat{e}_{ij}^q=\sum_{p=0}^q \hat{a}_{q,p}\hat{e}_{ji}^p$, we shall produce a recursion formula of $\{\hat{a}_{q,p}\}_{q\ge p\ge 0}$, where $\hat{a}_{q,p}$ is in the center $\mathfrak{Z}$. It follows from \eqref{eqn:3-16} that
\begin{equation}
\begin{split}
\hat{e}_{ij}^{q+1} &=\delta_{ji}\hat{c}_q-\hat{e}_{ji}^q-\sum_k \hat{e}_{kj}^q\hat{e}_{ki}\\
 &=\delta_{ji}\hat{c}_q-\hat{e}_{ji}^q-\sum_k \sum_p \hat{a}_{q,p}\hat{e}_{jk}^p\hat{e}_{ki} \\
&=\delta_{ji}\hat{c}_q-\hat{e}_{ji}^q-\sum_{p=0}^q \hat{a}_{q,p}\hat{e}_{ji}^{p+1}\\
 &=-\hat{a}_{q,q}\hat{e}_{ji}^{q+1}+(-\hat{a}_{q,q-1}-1)\hat{e}_{ji}^q-\sum_{p=0}^{q-2}\hat{a}_{q,p}\hat{e}_{ji}^{p+1}+\hat{c}_{q}\delta_{ji} \\
 &=\sum_{p=0}^{q+1}\hat{a}_{q+1,p}\hat{e}_{ji}^p. 
 \end{split}\nonumber
\end{equation}
Then we obtain a recursion formula for $\{\hat{a}_{q,p}\}_{q\ge p\ge 0}$, 
\begin{equation}
\hat{a}_{q+1,p}=
\begin{cases}
-\hat{a}_{q,q}  & p=q+1, \\
-\hat{a}_{q,q-1}-1  & p=q, \\
-\hat{a}_{q,p-1} & 1\le p\le q-1, \\
\hat{c}_q & p=0.
\end{cases}\nonumber
\end{equation}
Because of $\hat{e}_{ij}^0=\hat{e}_{ji}^0$ and $\hat{e}_{ij}^1=-\hat{e}_{ji}^1+(n-1)\hat{e}^0_{ij}$, the initial condition of $\hat{a}_{q,p}$ is $(\hat{a}_{0,0},\hat{a}_{1,0},\hat{a}_{1,1})=(1,n-1,-1)$. We solve the recursion formula and have 
\begin{equation}
\hat{a}_{q+1,p}=
\begin{cases}
(-1)^{q+1} & p=q+1, \\
(-1)^q(n-1)-\frac{1-(-1)^q}{2}=(-1)^q\hat{c}_0-\frac{1-(-1)^{q+1}}{2}  & p=q, \\(-1)^{p}\hat{c}_{q-p} & 0\le p\le q-1. 
\end{cases}\nonumber
\end{equation}
Thus we have proved the theorem. 
\end{proof}
In Section \ref{sec:6}, we show that  \eqref{eqn:3-17} induces all Bochner-Weitzenb\"ock formulas on Riemannian manifolds. Hence we call \eqref{eqn:3-17} \textit{the universal Bochner-Weitzenb\"ock formula}. 

Take the trace in \eqref{eqn:3-19}, and we have identities for $\{\hat{c}_q\}_{q\ge 0}$. 
\begin{cor}\label{cor:3-5}
The Casimir elements $\{\hat{c}_0,\hat{c}_1,\cdots \}$ satisfy 
\begin{equation}
2\hat{c}_{2q+1}=-\hat{c}_{2q}+\sum_{p=0}^{2q}(-1)^p\hat{c}_{2q-p}\hat{c}_p
\nonumber 
\end{equation}
for $q=0,1,\cdots$.
\end{cor}

D. Calderbank suggested to the author how the formula \eqref{eqn:3-17} is related with \cite{CGH}. We set 
\begin{equation}
\hat{E}_{ij}^q:=-\frac{1+(-1)^q}{2}\hat{e}_{ij}^q+\sum_{p=0}^q(-1)^p\hat{c}_{q-p}\hat{e}_{ij}^p. \label{sym}
\end{equation}
Then the universal Bochner-Weitzenb\"ock formula means $\hat{E}_{ij}^q=(-1)^q\hat{E}_{ji}^q$. This symmetry is the same as a formula in Theorem 4.8 of \cite{CGH}. The better point of our formula is that \eqref{eqn:3-17} is independent of representations, or universal. We notice that, if we find a formula with such a symmetry in the enveloping algebra, then we can obtain Bochner-Weitzenb\"ock formula. In fact, we have another universal Bochner-Weitzenb\"ock formula independent of \eqref{eqn:3-17} in the next section.

%%%%%%%%%%%%%%%%%%%%%%%%%%%%%%%%
%%            4               %%
%%%%%%%%%%%%%%%%%%%%%%%%%%%%%%%%
\section{Principal symbols of gradients}\label{sec:4}
We discuss the Clifford multiplication on spinor space $V_{\Delta}$, where $n$ is odd and $\Delta$ is $((1/2)_m)$. The Clifford multiplication is an action on $V_{\Delta}$ of $\xi$ in $\mathbb{R}^n$, 
\[ 
V_{\Delta}\ni \phi\mapsto \xi\cdot\phi\in V_{\Delta},
 \]
and satisfies the Clifford relation $\xi\cdot\eta\cdot+\eta\cdot\xi\cdot=-2\la\xi,\eta\ra \id$. To generalize the Clifford multiplication, we use another definition as follows. We consider the tensor representation $(\pi_{\Delta}\otimes \pi_{\mu_1},V_{\Delta}\otimes (\mathbb{R}^n\otimes \mathbb{C}))$ and decompose it into irreducible modules, $V_{\Delta}\otimes \mathbb{C}^n=V_T\oplus V_{\Delta}$. Here $T$ is $(3/2,(1/2)_{m-1})$. We denote by $\Pi_{\Delta}^{\Delta}$ the orthogonal projection from $V_{\Delta}\otimes \mathbb{C}^n$ onto $V_{\Delta}$ and show that $\Pi_{\Delta}^{\Delta} (\phi\otimes \xi)$ is a constant multiple of $\xi\cdot \phi$. Thus the projection $\Pi^{\Delta}_{\Delta}$ gives another definition  of the Clifford multiplication.

We generalize the above discussion to other representation spaces. We consider the tensor representation $(\pi_{\rho}\otimes\pi_{\mu_1},V_{\rho}\otimes \mathbb{C}^n)$ and its irreducible decomposition
\begin{equation}
V_{\rho}\otimes \mathbb{C}^n=\oplus_{\lambda} V_{\lambda}.
    \nonumber
\end{equation}
The highest weights of irreducible components occur with multiplicity one and are characterized as follows \cite{F}.
\begin{prop}
\begin{enumerate}
\item When $n=2m$, or when $n=2m+1$ and $\rho^m=0$, the highest weight of irreducible component in $V_{\rho}\otimes \mathbb{C}^n$ is dominant and $\rho\pm \mu_i$ for $i=1,2,\cdots,m$. 
\item When $n=2m+1$ and $\rho^m>0 $, the highest weight is dominant, and $\rho$ or $\rho \pm \mu_i$ for $i=1,2,\cdots,m$. 
\end{enumerate}
\end{prop}
\begin{ex}
\begin{enumerate}
	\item For $n=2m+1$ and $\rho=((1/2)_m)$, the highest weights in $V_{\rho}\otimes \mathbb{C}^n$ are $(3/2,(1/2)_{m-1})$ and $((1/2)_m)$. 
	\item For $n=2m$ and $\rho=((1/2)_{m-1},\pm 1/2)$, the highest weights in $V_{\rho}\otimes \mathbb{C}^n$ are $(3/2,(1/2)_{m-2},\pm 1/2)$ and $((1/2)_{m-1},\mp 1/2)$.
\end{enumerate}
\end{ex}
The tensor inner product on $V_{\rho}\otimes \mathbb{C}^n$ induces the one on $V_{\lambda}$. Hence each component is orthogonal to others. We denote by $\Pi_{\lambda}^{\rho}$ the orthogonal projection from $V_{\rho}\otimes \mathbb{C}^n$ onto $V_{\lambda}$. 
%%%%%%%%
\begin{defini}\label{def:4-1}
For $\xi$ in $\mathbb{C}^n$, we define a linear mapping $p^{\rho}_{\lambda}(\xi)$ from $V_{\rho}$ to $V_{\lambda}$ by
\begin{equation}
\mathbb{C}^n\times V_{\rho}\ni (\xi,\phi)\mapsto p^{\rho}_{\lambda}(\xi)\phi:=\Pi^{\rho}_{\lambda}(\phi\otimes \xi)\in V_{\lambda}. 
\nonumber 
\end{equation}
We denote by $p^{\rho}_{\lambda}(\xi)^{\ast}$ the adjoint operator of $p^{\rho}_{\lambda}(\xi)$ with respect to inner products on $V_{\rho}$ and $V_{\lambda}$. We call these linear mappings $p^{\rho}_{\lambda}(\xi)$ and $p^{\rho}_{\lambda}(\xi)^{\ast}$ \textit{the Clifford homomorphisms associated to $\rho$ and $\lambda$}.
\end{defini}
%%%%%%

We shall investigate properties of Clifford homomorphisms. 
%%%%%%%%%%%%%%%%%%%%%%%%%%%%%%%%
\begin{lem}\label{lem:4-2}
The Clifford homomorphism $p^{\rho}_{\lambda}$ satisfies 
\begin{equation}
\sum_i p^{\rho}_{\lambda}(e_i)\pi_{\rho}(e_{ij})=w(\rho;\lambda)p^{\rho}_{\lambda}(e_j) \label{eqn:4-3}
\end{equation}
for each $j$. Here, $w(\rho;\lambda)$ is a constant given by
\begin{equation}
w(\rho;\lambda):=1/2(\la \delta+\lambda,\delta +\lambda \ra-\la \delta+\rho,\delta+\rho\ra-n+1). 
\nonumber 
\end{equation}
We call this constant $w(\rho;\lambda)$ \textit{the conformal weight} associated to $\rho$ and $\lambda$.
\end{lem}
%%%%%%%%%%%%%%%%%%%%%%%%%%%%%%
\begin{proof}
We define an operator $C$ on $V_{\rho}\otimes \mathbb{C}^n$ by  
$$
C:=\pi_{\rho}\otimes
\pi_{\mu_1}(c_2)-\pi_{\rho}(c_2)\otimes \id -\id \otimes \pi_{\mu_1}(c_2).
$$
Since $\pi_{\rho}\otimes\pi_{\mu_1}(c_2)$ is $\pi_{\lambda}(c_2)$ on irreducible component $V_{\lambda}$, we show from \eqref{eqn:3-1} that $C$ is $4w(\rho;\lambda)\id$ on $V_{\lambda}$. Then we have 
$$C(\phi\otimes e_i)=C(\sum_{\lambda} p^{\rho}_{\lambda}(e_i)\phi)=\sum_{\lambda} 4w(\rho;\lambda)p^{\rho}_{\lambda}(e_i)\phi$$
for $\phi\otimes e_i$ in $V_{\rho}\otimes \mathbb{C}^n$. On the other hand, we know that 
\begin{equation}
\begin{split}
C=&\pi_{\rho}\otimes\pi_{\mu_1}(c_2)-\pi_{\rho}(c_2)\otimes \id -\id \otimes \pi_{\mu_1}(c_2) \\
=&\sum (\pi_{\rho}(e_{ij})\otimes \id+\id \otimes \pi_{\mu_1}(e_{ij}))(\pi_{\rho}(e_{ji})\otimes \id+\id \otimes \pi_{\mu_1}(e_{ji})) \\
&\quad -\sum \pi_{\rho}(e_{ij})\pi_{\rho}(e_{ji})\otimes \id -\sum \id \otimes \pi_{\mu_1}(e_{ij})\pi_{\mu_1}(e_{ji})\\
=&2\sum_{ij}\pi_{\rho}(e_{ij})\otimes \pi_{\mu_1}(e_{ji}).
\end{split}\nonumber
\end{equation}
Then, 
\begin{equation}
\begin{split}
C(\phi\otimes e_i)&=2\sum_{kl}\pi_{\rho}(e_{kl})\phi\otimes \pi_{\mu_1}(e_{lk})e_i =2\sum_{kl}\pi_{\rho}(e_{kl})\phi\otimes (\delta_{il}e_k-\delta_{ki}e_l)\\     &=4 \sum_{k} \pi_{\rho}(e_{ki})\phi\otimes e_k=4\sum_{\lambda}\sum_k p^{\rho}_{\lambda}(e_k)\pi_{\rho}(e_{ki})\phi.
\end{split}\nonumber
\end{equation}
As a result, we have $\sum_k p^{\rho}_{\lambda}(e_k)\pi_{\rho}(e_{ki})=w(\rho;\lambda)p^{\rho}_{\lambda}(e_i)$ for each $\lambda$. 
\end{proof}
%%%%%%%%%%%%%%%%%%%%%%%%%%%%%%%%%%%
From this lemma, we can relate Clifford homomorphisms to the enveloping algebra \cite{CGH}. 
%%%%%%%%%%%%%%%%%%%%%%%%%%%%%%%
\begin{prop}\label{prop:4-3}
The Clifford homomorphisms $\{p^{\rho}_{\lambda}\}_{\lambda}$ satisfy 
\begin{equation}
 \sum_{\lambda} w(\rho;\lambda)^{q}p^{\rho}_{\lambda}(e_i)^{\ast}p^{\rho}_{\lambda}(e_j)  =\pi_{\rho}(e_{ij}^{q}) \label{eqn:4-8}
\end{equation}
for $q=0,1,\cdots $, and $i,j=1,\cdots,n$. In particular, we have 
\begin{equation}
\sum_{\lambda} w(\rho;\lambda)^{q}\sum_i p^{\rho}_{\lambda}(e_i)^{\ast}p^{\rho}_{\lambda}(e_i)  =\pi_{\rho}(c_q). 
\nonumber 
\end{equation}
\end{prop}
%%%%%%%%%%
\begin{proof}
For $\phi$ and $\psi$ in $V_{\rho}$, we have 
$$
\delta_{ij}\la \phi,\psi\ra=\la \phi\otimes e_i,\psi\otimes e_j\ra=\sum_{\lambda}\la p^{\rho}_{\lambda}(e_i)\phi,p^{\rho}_{\lambda}(e_j)\psi\ra=\la \sum_{\lambda}p^{\rho}_{\lambda}(e_j)^{\ast}p^{\rho}_{\lambda}(e_i)\phi,\psi\ra .
$$
Then we have proved \eqref{eqn:4-8} for $q=0$, 
\begin{equation}
\sum_{\lambda}p^{\rho}_{\lambda}(e_j)^{\ast}p^{\rho}_{\lambda}(e_i)=\delta_{ji}.   \label{eqn:4-10}
\end{equation}
By using this equation and \eqref{eqn:4-3}, we can prove \eqref{eqn:4-8} inductively. 
\end{proof}
%%%%%%%%%%%%%%%%
We calculate the conformal weights for $\lambda=\rho,\rho\pm \mu_i$, 
\begin{equation}
\begin{cases}
w(\rho;\rho+\mu_i)=\rho^i+1-i,   &i=1,\cdots,m, \\
w(\rho;\rho-\mu_i)=-\rho^i-n+i+1, &i=1,\cdots,m, \\
w(\rho;\rho)=-(n-1)/2.           
\end{cases}\label{eqn:4-5}
\end{equation}
We assume that the number of irreducible components is $N$, that is, $N=\#\{\lambda|V_{\lambda}\subset V_{\rho}\otimes \mathbb{C}^n\}$. Arranging them lexicographically  as $\lambda_1=\rho+\mu_1>\lambda_2>\cdots>\lambda_N$, we show from \eqref{eqn:4-5} that
\begin{equation}
w(\rho;\lambda_1)>w(\rho;\lambda_2)>\cdots >w(\rho;\lambda_N) 
\nonumber 
\end{equation}
except the following case. When $n=2m$, $\rho^{m-1}>0$ and $\rho^m=0$, there always exist the highest weights $\lambda_+:=\rho+\mu_m$ and $\lambda_-:=\rho-\mu_m$, whose conformal weights coincide. Then we have 
\begin{equation}
w(\rho;\lambda_1)>w(\rho;\lambda_2)>\cdots> w(\rho;\lambda_+)=w(\rho;\lambda_-)>\cdots >w(\rho;\lambda_N). \nonumber 
\end{equation}
We call this case \textit{the exceptional case}. Thus the conformal weights differ from each other except the exceptional case. It follows from Proposition \ref{prop:4-3} that $p^{\rho}_{\lambda}(e_i)^{\ast}p^{\rho}_{\lambda}(e_j)$ can be realized as a linear combination of $\{\pi_{\rho}(e_{ij}^q)\}_q$. In fact, ordering highest weights as $\lambda_1>\cdots >\lambda_N$, we obtain 
\begin{equation}
(p^{\rho}_{\lambda_1}(e_i)^{\ast}p^{\rho}_{\lambda_1}(e_j),\cdots, p^{\rho}_{\lambda_N}(e_i)^{\ast}p^{\rho}_{\lambda_N}(e_j))W^t=(\delta_{ij},\pi_{\rho}(e_{ij}),\cdots, \pi_{\rho}(e_{ij}^{N-1})), \nonumber
\end{equation}
where $W$ is a $N\times N$ Vandermonde matrix, 
\begin{equation}
\begin{pmatrix}
1 &1& \cdots & 1\\
w(\rho;\lambda_1)& w(\rho;\lambda_2) &\cdots & w(\rho;\lambda_N) \\
\cdots & \cdots &\cdots& \cdots \\
w(\rho;\lambda_1)^{N-1}& w(\rho;\lambda_2)^{N-1} &\cdots & w(\rho;\lambda_N)^{N-1} 
\end{pmatrix}. \nonumber
\end{equation}
Since $W$ is invertible except the exceptional case, the Clifford homomorphism $p^{\rho}_{\lambda}(e_i)^{\ast}p^{\rho}_{\lambda}(e_j)$ is a linear combination of $\{\pi_{\rho}(e_{ij}^q)\}_q$. 

For the exceptional case, we need the Pfaffian element $\pf$ in \eqref{eqn:3-7} to distinguish $p^{\rho}_{\lambda_+}$ from $p^{\rho}_{\lambda_-}$. So we shall investigate relations between the Pfaffian element and Clifford homomorphisms. 
%%%%%%%%%%%%
\begin{defini}
We define an element $\pf_{ij}$ in the enveloping algebra for $i,j=1,\cdots, 2m$ by 
\begin{equation}
\pf_{ij}:=\left\{
\begin{array}{lc}
 \pf,  & i=j, \\
 (-1)^{i+j}\frac{2m}{(\sqrt{-1})^m2^mm!} \sum_{\sigma\in \mathfrak{S}^{ij}_{2m}} \mathrm{sgn}(\sigma)e_{\sigma(1)\sigma(2)}\cdots e_{\sigma(2m-1)\sigma(2m)}, &i<j, \\
-\pf_{ji} &i>j,
\end{array}\right.
\nonumber 
\end{equation}
where $\mathfrak{S}^{ij}_{2m}$ is the permutation group of $\{1,\cdots,2m\} \setminus \{i,j\}$. 
\end{defini}
\begin{ex}[four dimensional case]\label{ex:4-1}
Let $\{e_i\}_{i=1}^4$ be an oriented orthonormal basis of $\mathbb{R}^4$. We calculate $\pf_{ij}$,
\begin{equation}
\begin{split}
\pf_{12}=e_{34}, \quad \pf_{13}=-e_{24}, \quad \pf_{14}=e_{23}, \\
\pf_{23}=e_{14}, \quad \pf_{24}=-e_{13}, \quad \pf_{34}=e_{12}.
\end{split}\label{eqn:pf-4}
\end{equation}
We denote the Hodge operator by $\ast:\Lambda^2(\mathbb{R}^4)\to \Lambda^2(\mathbb{R}^4)$. Then, $\pf_{ij}=\ast e_{ij}$ for $i\neq j$. 
\end{ex}
From the definition of $\pf_{ij}$, we have the following. 
\begin{prop}
The elements $\{\pf_{ij}\}_{i,j}$ satisfy 
\begin{equation}
\pf_{ij}+\pf_{ji}=2\delta_{ij}\pf, \label{eqn:4-13} 
\end{equation}
and $\sum_i\pf_{ii}=2m \pf$.
\end{prop}
The equation \eqref{eqn:4-13} induces a symmetry $\pf_{ij}-\delta_{ij}\pf=-(\pf_{ji}-\delta_{ji}\pf)$ as \eqref{sym} and gives Bochner-Weitzenb\"ock formula in Section \ref{sec:6}. So we call \eqref{eqn:4-13} the universal Bochner-Weitzenb\"ock formula as well as \eqref{eqn:3-17}.
%%%%%%%%%%%%%%%
%%%%%%%%%
\begin{prop}\label{prop:4-7}
We can connect $\pf_{ij}$ to the Clifford homomorphisms $\{p^{\rho}_{\lambda}\}_{\lambda}$ as 
\begin{equation}
\sum_{\lambda} \pi_{\lambda}(\pf)p^{\rho}_{\lambda}(e_i)^{\ast}p^{\rho}_{\lambda}(e_j)=\pi_{\rho}(\pf_{ij}).
\label{eqn:4-14}
\end{equation}
In particular, we have
\begin{equation}
\sum_{\lambda}\pi_{\lambda}(\pf)\sum_i p^{\rho}_{\lambda}(e_i)^{\ast}p^{\rho}_{\lambda}(e_i)=2m\pi_{\rho}(\pf). \label{eqn:4-15}
\end{equation}
\end{prop}
%%%%%%%%
%%%%%%%%
\begin{proof}
For $\phi$ and $\psi$ in $V_{\rho}$, 
\begin{equation}
\begin{split}
&\la\pi_{\rho}\otimes \pi_{\mu_1}(\pf)(\phi\otimes e_j),\psi \otimes e_i\ra =\sum_{\lambda}\la\pi_{\lambda}(\pf)p^{\rho}_{\lambda}(e_j)\phi, p^{\rho}_{\lambda}(e_i)\psi \ra \\
=&\la\sum_{\lambda}\pi_{\lambda}(\pf)p^{\rho}_{\lambda}(e_i)^{\ast}p^{\rho}_{\lambda}(e_j)\phi, \psi \ra.
\end{split}\label{eqn:4-16}
\end{equation}
On the other hand, we know 
\begin{equation}
\begin{split}
&\pi_{\rho}\otimes \pi_{\mu_1}(e_{\sigma(2l-1)\sigma(2l)})(\phi\otimes e_j)\\
=&
\begin{cases}
\pi_{\rho}(e_{\sigma(2l-1)\sigma(2l)})\phi\otimes e_j+\phi\otimes e_{\sigma(2l)} , & \sigma(2l-1)=j, \\
\pi_{\rho}(e_{\sigma(2l-1)\sigma(2l)})\phi\otimes e_j-\phi\otimes e_{\sigma(2l-1)}, & \sigma (2l)=j, \\
\pi_{\rho}(e_{\sigma(2l-1)\sigma(2l)})\phi\otimes e_j ,& \textrm{otherwise}.
\end{cases}
\end{split}\nonumber
\end{equation}
Then taking account of $\la e_k,e_l\ra=\delta_{kl}$, we have
\begin{equation}
\begin{split}
&(\sqrt{-1})^m2^mm! \la\pi_{\rho}\otimes \pi_{\mu_1}(\pf)(\phi\otimes e_j),\psi \otimes e_i\ra \\
=&\la\pi_{\rho}\otimes \pi_{\mu_1}(\sum_{\sigma\in \mathfrak{S}_{2m}} \mathrm{sgn}(\sigma)e_{\sigma(1)\sigma(2)}\cdots e_{\sigma(2m-1)\sigma(2m)})(\phi\otimes e_j),\psi\otimes e_i\ra \\=&
\begin{cases}
\la\pi_{\rho}((\sqrt{-1})^m2^mm! \pf) \phi,\psi\ra, & i=j,\\
(-1)^{i+j}2m\la\pi_{\rho}(\sum_{\sigma\in \mathfrak{S}^{ij}_{2m}} \mathrm{sgn}(\sigma)e_{\sigma(1)\sigma(2)}\cdots e_{\sigma(2m-1)\sigma(2m)})\phi,\psi\ra, &i<j,\\
(-1)^{i+j+1}2m\la\pi_{\rho}(\sum_{\sigma\in \mathfrak{S}^{ij}_{2m}} \mathrm{sgn}(\sigma)e_{\sigma(1)\sigma(2)}\cdots e_{\sigma(2m-1)\sigma(2m)})\phi,\psi\ra, &i>j.
\end{cases}
\end{split}\nonumber 
\end{equation}
Combining this equation and \eqref{eqn:4-16}, we have \eqref{eqn:4-14}. The equation \eqref{eqn:4-15} is easy from \eqref{eqn:4-13}.
\end{proof}
%%%%%%%%%
We return to  the exceptional case. It follows from \eqref{eqn:3-8} that $\pi_{\lambda_+}(\pf)=-\pi_{\lambda_-}(\pf)\neq 0$ and $\pi_{\lambda}(\pf)=0$ for $\lambda\neq \lambda_{\pm}$. Then
\begin{equation}
p^{\rho}_{\lambda_+}(e_i)^{\ast}p^{\rho}_{\lambda_+}(e_j)-p^{\rho}_{\lambda_-}(e_i)^{\ast}p^{\rho}_{\lambda_-}(e_j)=\frac{1}{\pi_{\lambda_+}(\pf)}\pi_{\rho}(\pf_{ij}).
 \nonumber 
\end{equation}
Therefore the Clifford homomorphism $p^{\rho}_{\lambda_{\pm}}(e_i)^{\ast}p^{\rho}_{\lambda_{\pm}}(e_j)$ is a linear combination of $\{\pi_{\rho}(e_{ij}^q)\}_q\cup \{\pi_{\rho}(\pf_{ij})\}$. 

As a result,  we have the following corollary of Proposition \ref{prop:4-3} and \ref{prop:4-7}. 
%%%%%%%%%%%%%
\begin{cor}
We can rewrite  $p^{\rho}_{\lambda}(e_i)^{\ast}p^{\rho}_{\lambda}(e_j)$ as a linear combinations of $\{\pi_{\rho}(e_{ij}^q)\}_q$ and $\pi_{\rho}(\pf_{ij})$. 
\end{cor}
%%%%%%%%%
This corollary implies that we may investigate $e_{ij}^q$ and $\pf_{ij}$ instead of  Clifford homomorphisms.  In fact, the universal Bochner-Weitzenb\"ock formulas \eqref{eqn:3-17} and \eqref{eqn:4-13} give a lot of relations for Clifford homomorphisms. 

First, we consider the equation \eqref{eqn:4-8} for $q=0,1$. Taking account of $e_{ij}=-e_{ji}$,  we have 
\begin{gather}
\sum_{\lambda}(p^{\rho}_{\lambda}(e_i)^{\ast}p^{\rho}_{\lambda}(e_j)+p^{\rho}_{\lambda}(e_j)^{\ast}p^{\rho}_{\lambda}(e_i))=2\delta_{ij},  \label{eqn:4-20}
\\
\sum_{\lambda}w(\rho;\lambda)(p^{\rho}_{\lambda}(e_i)^{\ast}p^{\rho}_{\lambda}(e_j)+p^{\rho}_{\lambda}(e_j)^{\ast}p^{\rho}_{\lambda}(e_i))=0. \label{eqn:4-21}
\end{gather}
Note that we can obtain the Clifford relation on spinor space from the above two equations.  

To construct further relations, we use the universal Bochner-Weitzenb\"ock formulas \eqref{eqn:3-17} and \eqref{eqn:4-13}. We define the translated conformal weight $\hat{w}(\rho;\lambda)$ by 
\begin{equation}
\hat{w}(\rho;\lambda):=w(\rho;\lambda)+\frac{n-1}{2}. 
\nonumber 
\end{equation}
Then it is easy shown that 
\begin{equation}
\sum_{\lambda}\hat{w}(\rho;\lambda)^qp^{\rho}_{\lambda}(e_i)^{\ast}p^{\rho}_{\lambda}(e_j)=\pi_{\rho}(\hat{e}_{ij}^q). \nonumber 
\end{equation}
Substituting \eqref{eqn:3-18} for this equation, we have 
\begin{equation}
\begin{split}
&\sum_{\lambda}\hat{w}(\rho;\lambda)^{2q}p^{\rho}_{\lambda}(e_i)^{\ast}p^{\rho}_{\lambda}(e_j)\\
=&\sum_{\lambda}\{\hat{w}(\rho;\lambda)^{2q}+\sum_{p=0}^{2q-1}(-\hat{w}(\rho;\lambda))^p\pi_{\rho}(\hat{c}_{2q-1-p})\}p^{\rho}_{\lambda}(e_j)^{\ast}p^{\rho}_{\lambda}(e_i),
\end{split}\nonumber
\end{equation}
and hence, 
\begin{gather}
\sum_{\lambda}\{\sum_{p=0}^{2q-1}(-\hat{w}(\rho;\lambda))^p\pi_{\rho}(\hat{c}_{2q-1-p})\}(p^{\rho}_{\lambda}(e_i)^{\ast}p^{\rho}_{\lambda}(e_j)+p^{\rho}_{\lambda}(e_j)^{\ast}p^{\rho}_{\lambda}(e_i))=0. \nonumber 
\end{gather}

From \eqref{eqn:4-13} and \eqref{eqn:4-14} for $n=2m$, we have 
\begin{equation}
\sum_{\lambda}\pi_{\lambda}(\pf)(p^{\rho}_{\lambda}(e_i)^{\ast}p^{\rho}_{\lambda}(e_j)+p^{\rho}_{\lambda}(e_j)^{\ast}p^{\rho}_{\lambda}(e_i))=2\pi_{\rho}(\pf)\delta_{ij}. \nonumber 
\end{equation}
We can easily show from \eqref{eqn:3-8} and \eqref{eqn:4-5} that 
\begin{equation}
(w(\rho;\lambda)+m-1)\pi_{\lambda}(\pf)=(w(\rho;\lambda)+m)\pi_{\rho}(\pf) \nonumber
\end{equation}
for each $\lambda$. Then 
\begin{equation}
\begin{split}
 &\sum_{\lambda}\pi_{\lambda}(\pf)w(\rho;\lambda)(p^{\rho}_{\lambda}(e_i)^{\ast}p^{\rho}_{\lambda}(e_j)+p^{\rho}_{\lambda}(e_j)^{\ast}p^{\rho}_{\lambda}(e_i))\\
=&\sum_{\lambda}(-(m-1)\pi_{\lambda}(\pf)+w(\rho;\lambda)\pi_{\rho}(\pf)+m\pi_{\rho}(\pf))(p^{\rho}_{\lambda}(e_i)^{\ast}p^{\rho}_{\lambda}(e_j)+p^{\rho}_{\lambda}(e_j)^{\ast}p^{\rho}_{\lambda}(e_i))\\
=&-2(m-1)\pi_{\rho}(\pf)\delta_{ij}+2m\pi_{\rho}(\pf)\delta_{ij}=2\pi_{\rho}(\pf)\delta_{ij}.
\end{split}\nonumber
\end{equation}
As a result, we have
\begin{equation}
\sum_{\lambda}\pi_{\lambda}(\pf)(w(\rho;\lambda)-1)(p^{\rho}_{\lambda}(e_i)^{\ast}p^{\rho}_{\lambda}(e_j)+p^{\rho}_{\lambda}(e_j)^{\ast}p^{\rho}_{\lambda}(e_i))=0. \nonumber 
\end{equation}
Thus we obtain relations among Clifford homomorphisms like the Clifford relation as follows. 
%%%%%%%%%%%%%%%%%%%%%%%%%%%%%
\begin{thm}\label{thm:4-7}
Let $\xi$ and $\eta$ be in $\mathbb{R}^n$. Then the Clifford homomorphisms $\{p^{\rho}_{\lambda}\}_{\lambda}$ satisfy that, 
\begin{gather}
\sum_{\lambda}(p^{\rho}_{\lambda}(\xi)^{\ast}p^{\rho}_{\lambda}(\eta)+p^{\rho}_{\lambda}(\eta)^{\ast}p^{\rho}_{\lambda}(\xi))=2\la\xi,\eta\ra, \nonumber \\
\sum_{\lambda}\{\sum_{p=0}^{2q-1}(-\hat{w}(\rho;\lambda))^p\pi_{\rho}(\hat{c}_{2q-1-p})\}(p^{\rho}_{\lambda}(\xi)^{\ast}p^{\rho}_{\lambda}(\eta)+p^{\rho}_{\lambda}(\eta)^{\ast}p^{\rho}_{\lambda}(\xi))=0, \nonumber 
\end{gather}
and, for $n=2m$, 
\begin{gather}
\sum_{\lambda}\pi_{\lambda}(\pf)(p^{\rho}_{\lambda}(\xi)^{\ast}p^{\rho}_{\lambda}(\eta)+p^{\rho}_{\lambda}(\eta)^{\ast}p^{\rho}_{\lambda}(\xi))=2\pi_{\rho}(\pf)\la\xi,\eta\ra, \label{eqn:4-29}\\
\sum_{\lambda}\pi_{\lambda}(\pf)(w(\rho;\lambda)-1)(p^{\rho}_{\lambda}(\xi)^{\ast}p^{\rho}_{\lambda}(\eta)+p^{\rho}_{\lambda}(\eta)^{\ast}p^{\rho}_{\lambda}(\xi))=0. \label{eqn:4-30}
\end{gather}
\end{thm}
%%%%%%%%%%%%%%%%%%%%%%%%%%%%%%

In the rest of this section, we calculate the eigenvalues of $c_q$ on irreducible $\mathfrak{so}(n)$-modules. Our method is based on \cite{CGH} and \cite{O}. It is known that the usual Clifford multiplication satisfies $(\pi_{\mu_1}(g)\xi)\cdot=\pi_{\Delta}(g)\xi \cdot \pi_{\Delta}(g^{-1})$ for $g$ in $Spin(n)$ and $\xi$ in $\mathbb{R}^n$. Clifford homomorphism is also compatible with the action of $SO(n)$ or $Spin(n)$. 
\begin{lem}\label{lem:4-8}
For $g$ in $SO(n)$ or $Spin(n)$ and $\xi$ in $\mathbb{C}^n$, we have 
\begin{equation}
p^{\rho}_{\lambda}(\pi_{\mu_1}(g)\xi)=\pi_{\lambda}(g)p^{\rho}_{\lambda}(\xi)\pi_{\rho}(g^{-1}). \label{eqn:4-31}
\end{equation}
Hence, for $e_{ij}$ in $\mathfrak{so}(n)$,  
\begin{equation}
p^{\rho}_{\lambda}(\pi_{\mu_1}(e_{ij})\xi)=
\pi_{\lambda}(e_{ij})p^{\rho}_{\lambda}(\xi)-p^{\rho}_{\lambda}(\xi)\pi_{\rho}(e_{ij}). \label{eqn:4-32}
\end{equation}
\end{lem}
\begin{proof}
Consider the action of $g$ on $V_{\rho}\otimes \mathbb{C}^n=\oplus_{\lambda}V_{\lambda}$, and we have 
\begin{equation}
\sum_{\lambda} \pi_{\lambda}(g)p^{\rho}_{\lambda}(\xi)\phi=\pi_{\rho}\otimes \pi_{\mu_1}(g)(\phi\otimes \xi)=\pi_{\rho}(g)\otimes \pi_{\mu_1}(g)\xi=\sum_{\lambda}p^{\rho}_{\lambda}(\pi_{\mu_1}(g)\xi)\pi_{\rho}(g)\phi\nonumber
\end{equation}
for $\phi\otimes \xi$ in $V_{\rho}\otimes \mathbb{C}^n$. Then we conclude $p^{\rho}_{\lambda}(\pi_{\mu_1}(g)\xi)=\pi_{\lambda}(g)p^{\rho}_{\lambda}(\xi)\pi_{\rho}(g^{-1})$. Its infinitesimal realization is \eqref{eqn:4-32}. 
\end{proof}
The Clifford homomorphism $p^{\rho}_{\lambda}$ is defined through the projection $\Pi^{\rho}_{\lambda}:V_{\rho}\otimes \mathbb{C}^n\to V_{\lambda}$. Therefore $\Pi^{\rho}_{\lambda}$ is realized with the Clifford homomorphism. 
%%%%%%%%%%%
\begin{lem}\label{lem:4-9}
The orthogonal projection $\Pi^{\rho}_{\lambda}:V_{\rho}\otimes \mathbb{C}^n\to V_{\lambda}\subset V_{\rho}\otimes \mathbb{C}^n$ is realized as follows.
\begin{equation}
\Pi^{\rho}_{\lambda}(\phi\otimes \xi)=\sum_i p^{\rho}_{\lambda}(e_i)^{\ast}p^{\rho}_{\lambda}(\xi)\phi\otimes e_i. \label{eqn:4-33}
\end{equation}
\end{lem}
%%%%%%%%%%%
\begin{proof}
It follows from  \eqref{eqn:4-31} that the following mapping is an $\mathfrak{so}(n)$-equivariant injection,
$$
V_{\lambda}\ni \psi\mapsto \sum_i p^{\rho}_{\lambda}(e_i)^{\ast}\psi\otimes e_i \in V_{\rho}\otimes \mathbb{C}^n.
$$
Taking account of \eqref{eqn:4-10}, we decompose $\phi\otimes \xi$,
\begin{equation}
\begin{split}
\phi\otimes \xi &=\sum_{i}\langle \xi, e_i\rangle\phi\otimes e_i=\sum_{i}\sum_{\lambda} p^{\rho}_{\lambda}(e_i)^{\ast}p^{\rho}_{\lambda}(\xi)\phi\otimes e_i \\
    &=\sum_{\lambda}\sum_{i} p^{\rho}_{\lambda}(e_i)^{\ast}p^{\rho}_{\lambda}(\xi)\phi\otimes e_i.
\end{split}\nonumber
\end{equation}
Since $\sum_{i} p^{\rho}_{\lambda}(e_i)^{\ast}p^{\rho}_{\lambda}(\xi)\phi\otimes e_i$ is in $V_{\lambda}$ for each $\lambda$, we find out the projection formula \eqref{eqn:4-33}.
\end{proof}
Lemma \ref{lem:4-8} leads that $\sum_i p^{\rho}_{\lambda}(e_i)^{\ast}p^{\rho}_{\lambda}(e_i)$ is compatible with the action of $\mathfrak{so}(n)$ and constant on $V_{\rho}$. 
%%%%%%%%%%%%%%%%%%%%
\begin{prop}\label{prop:4-10}
We set $d(\rho):=\dim V_{\rho}$. Then 
\begin{equation}
\sum_i p^{\rho}_{\lambda}(e_i)^{\ast}p^{\rho}_{\lambda}(e_i)=d(\lambda)/d(\rho). \nonumber 
\end{equation}
The eigenvalues of $c_q$ and $\hat{c}_q$ on irreducible $\mathfrak{so}(n)$-module $V_{\rho}$ are
\begin{equation}
\pi_{\rho}(c_q)=\frac{1}{d(\rho)}\sum_{\lambda}w(\rho;\lambda)^q d(\lambda), \quad \pi_{\rho}(\hat{c}_q)=\frac{1}{d(\rho)}\sum_{\lambda}\hat{w}(\rho;\lambda)^q d(\lambda). \label{eqn:4-35}
\end{equation}
Moreover, we have a relation for eigenvalues of the Pfaffian element,
\begin{equation}
2m\pi_{\rho}(\pf)=\frac{1}{d(\rho)}\sum_{\lambda}\pi_{\lambda}(\pf)d(\lambda) =\frac{1}{d(\rho)}\sum_{\lambda}\pi_{\lambda}(\pf)w(\rho;\lambda)d(\lambda). 
\label{eqn:4-36}
\end{equation}
\end{prop}
%%%%%%%%%%%%%%%
\begin{proof}
Let $\{\phi_{\alpha}\}_{\alpha=1}^{\dim V_{\rho}}$ be an orthonormal basis of $V_{\rho}$. Taking the trace of $\Pi^{\rho}_{\lambda}$,  we have
\begin{equation}
\begin{split}
 d(\lambda)=&\sum_{\alpha,i} \la \Pi^{\rho}_{\lambda}(\phi_{\alpha}\otimes e_i), \phi_{\alpha}\otimes e_i\ra =\sum_{\alpha,i} \la\sum_j p_{\lambda}^{\rho}(e_j)^{\ast}p^{\rho}_{\lambda}(e_i)(\phi_{\alpha}) \otimes e_j ,\phi_{\alpha}\otimes e_i\ra \\ 
=&\sum_{\alpha,i} \la \sum_j p_{\lambda}^{\rho}(e_j)^{\ast}p^{\rho}_{\lambda}(e_i)(\phi_{\alpha}),\phi_{\alpha}\ra\delta_{ij} =\sum_{\alpha} \la \phi_{\alpha},\phi_{\alpha}\ra\sum_i p_{\lambda}^{\rho}(e_i)^{\ast}p^{\rho}_{\lambda}(e_i) \\=&d(\rho) \sum_i p_{\lambda}^{\rho}(e_i)^{\ast}p^{\rho}_{\lambda}(e_i). \end{split}\nonumber
\end{equation}
Thus we obtain $\sum_i p^{\rho}_{\lambda}(e_i)^{\ast}p^{\rho}_{\lambda}(e_i)=d(\lambda)/d(\rho)$ and easily show \eqref{eqn:4-35} and \eqref{eqn:4-36}.
\end{proof}
\begin{rem}
The dimension $d(\rho)$ is calculated by Weyl's dimension formula \cite{K}, \cite{Z}. A variety of formulas for $\pi_{\rho}(c_q)$ in \cite{CGH}, \cite{NR} and \cite{O} are useful for explicit calculations.  
\end{rem}

%%%%%%%%%%%%%%%%%%%%%%%%%%%%%%%%
%%            5               %%
%%%%%%%%%%%%%%%%%%%%%%%%%%%%%%%%
\section{Gradients on Riemannian manifolds}\label{sec:5}
In this section, we define gradients and study their fundamental properties. We consider only gradients on Riemannian manifolds. The spin case is left to the readers, where the spin connection is used instead of the Levi-Civita connection \cite{Fr}, \cite{LM}. 

Let $(M,g)$ be an $n$-dimensional oriented Riemannian manifold, and $\textbf{SO}(M)$ be the principal $SO(n)$ bundle of the oriented orthonormal frames on $M$. For an irreducible unitary representation $(\pi_{\rho},V_{\rho})$ of $SO(n)$,  we have an associated Hermitian vector bundle $\mathbf{S}_{\rho}:=\mathbf{SO}(M)\times_{\pi_{\rho}}V_{\rho}$. The Levi-Civita connection on $\mathbf{SO}(M)$ gives a covariant derivative $\nabla$ on $\mathbf{S}_{\rho}$ compatible with fiber metric as follows. Let $e=(e_1,\cdots,e_n)$ be a local section of $\mathbf{SO}(M)$. For a unitary basis $\{\phi_{\alpha}\}_{\alpha}$ of $V_{\rho}$,  we have a local frame $\{[e,\phi_{\alpha}]\}_{\alpha}$ of $\mathbf{S}_{\rho}$. With respect to this local trivialization, a covariant derivative $\nabla$ is defined to be 
\begin{equation}
\nabla:=d+\frac{1}{2}\sum_{ij}g(\nabla^T e_i, e_j)\pi_{\rho}(e_{ij}), \label{eqn:5-1}
\end{equation}
where $\nabla^T$ is the Levi-Civita connection on the tangent bundle $T(M)$. Since the connection $1$-from is skew Hermitian, the derivative $\nabla$ is compatible with fiber metric, that is, $X\la\phi,\psi\ra=\la\nabla_X\phi,\psi\ra+\la\phi,\nabla_X\psi\ra$ for every vector field $X$.

We shall extend Clifford homomorphisms to bundle homomorphisms. We consider the tensor bundle $\mathbf{S}_{\rho}\otimes T_{\mathbb{C}}(M)=\mathbf{S}_{\rho}\otimes (T(M)\otimes\mathbb{C})$ and decompose it as
\begin{equation}
\mathbf{S}_{\rho}\otimes T_{\mathbb{C}}(M)=\oplus_{\lambda} \mathbf{S}_{\lambda}. \label{eqn:5-10}
\end{equation}
For each vector field $X=\sum X^i e_i$, we define a bundle homomorphism $p^{\rho}_{\lambda}(X)$ in $\Gamma (M,\Hom (\mathbf{S}_{\rho},\mathbf{S}_{\lambda}))$ by \begin{equation}
p^{\rho}_{\lambda}(X):\mathbf{S}_{\rho}\ni [e,\phi]\mapsto \sum_i X^i[e, p^{\rho}_{\lambda}(e_i)\phi]\in \mathbf{S}_{\lambda}.\nonumber 
\end{equation}
From \eqref{eqn:4-31}, we know that this bundle homomorphism is well-defined.  Furthermore, we can show from \eqref{eqn:4-32} and \eqref{eqn:5-1} that 
\begin{equation}
\nabla_Y(p^{\rho}_{\lambda}(X)\phi)=p^{\rho}_{\lambda}(\nabla^T_X Y)\phi+p^{\rho}_{\lambda}(X)\nabla_Y\phi \nonumber 
\end{equation}
for $\phi$ in $\Gamma(M,\mathbf{S}_{\rho})$. 

We define geometric first order differential operators depending on the metric $g$ on each associated bundle. 
\begin{defini}
We decompose $\nabla$ along \eqref{eqn:5-10}. Then we have a first order differential operator $D^{\rho}_{\lambda}=\Pi^{\rho}_{\lambda}\circ \nabla$, 
\begin{equation}
D^{\rho}_{\lambda}:\Gamma (M,\mathbf{S}_{\rho})\xrightarrow{\nabla}\Gamma (M,\mathbf{S}_{\rho}\otimes T^{\ast}_{\mathbb{C}}(M))\xrightarrow{\simeq}\Gamma (M,\mathbf{S}_{\rho}\otimes T_{\mathbb{C}}(M))\xrightarrow{\Pi^{\rho}_{\lambda}}\Gamma (M,\mathbf{S}_{\lambda}) \nonumber 
\end{equation}
for each $\lambda$. Here $\Pi^{\rho}_{\lambda}$ is the orthogonal projection defined fiberwise from $\mathbf{S}_{\rho}\otimes T_{\mathbb{C}}(M)$ onto $\mathbf{S}_{\lambda}$. We call this first order differential operator $D^{\rho}_{\lambda}$ \textit{the gradient associated to $\rho$ and $\lambda$}. 
\end{defini}
\begin{ex}
Let $\mathbf{S}_{\Delta}$ be the spinor bundle, where $n=2m+1$ and $\Delta=((1/2)_m)$. We have the irreducible decomposition $\mathbf{S}_{\Delta}\otimes T_{\mathbb{C}}(M)=\mathbf{S}_{T}\oplus \mathbf{S}_{\Delta}$ with $T=(3/2,(1/2)_{m-1})$. Then $D^{\Delta}_{\Delta}$ is the Dirac operator and $D^{\Delta}_T$ is the twistor operator up to a normalization. 
\end{ex}
\begin{ex}
Let $\Lambda^p(M)\otimes \mathbb{C}=\mathbf{S}_{(1_p)}$ be  the bundle of differential forms for $1\le p\le [n/2]$. The irreducible decomposition of $\Lambda^p(M)\otimes T_{\mathbb{C}}(M)$ is $\mathbf{S}_{(2,1_{p-1})}\oplus \mathbf{S}_{(1_{p+1})}\oplus \mathbf{S}_{(1_{p-1})}$. Then we have the conformal Killing operator $C$, the exterior derivative $d$, and the interior derivative $d^{\ast}$ up to a normalization.
\end{ex}
Because the principal symbol of $D^{\rho}_{\lambda}$ is the Clifford homomorphism $p^{\rho}_{\lambda}$, we find a formula of the gradient $D^{\rho}_{\lambda}$,
\begin{equation}
D^{\rho}_{\lambda}(\phi)=\Pi^{\rho}_{\lambda}(\sum_i \nabla_{e_i}\phi\otimes e_i^{\ast})=\sum_i p^{\rho}_{\lambda}(e_i)\nabla_{e_i}\phi.
\label{eqn:5-14}
\end{equation}
From a similar discussion to the Dirac operator \cite{LM}, we show that the formal adjoint operator $(D^{\rho}_{\lambda})^{\ast}$ of $D^{\rho}_{\lambda}$ is 
\begin{equation}
(D^{\rho}_{\lambda})^{\ast}=-\sum_i p^{\rho}_{\lambda}(e_i)^{\ast}\nabla_{e_i}.\nonumber 
\end{equation}

An important feature of the Dirac operator is conformal covariance. The gradients are also conformally covariant operators. Though this fact has been shown by H. D. Fegan \cite{F}, we give an explicit proof. We deform the Riemannian metric $g$ conformally as $g'=\exp (2\sigma)g$ for $\sigma$ in $C^{\infty}(M)$. We denote the objects associated to $g'$ by adding a symbol $``\; '\; "$ to them. The orthonormal frame bundle $\mathbf{SO}(M)$ is isomorphic to $\mathbf{SO}'(M)$ as a principal bundle, 
   \begin{equation}
   \Phi:\mathbf{SO}(M)\ni e=(e_1,\cdots,e_n)\mapsto e'=e^{-\sigma}(e_1,\cdots,e_n)  \in \mathbf{SO}'(M). \nonumber 
    \end{equation}
Then there is a bundle isometry for each $\rho$, 
  \begin{equation}
  \Phi_{\rho}:\mathbf{S}_{\rho}=\mathbf{SO}(M)\times_{\pi_{\rho}}V_{\rho}\ni [e,\phi]\mapsto [e',\phi]\in \mathbf{SO}'(M)\times_{\pi_{\rho}}V_{\rho}=\mathbf{S}'_{\rho} \nonumber
  \end{equation}
 such that $\Phi_{\lambda}\circ p^{\rho}_{\lambda}(X)=e^{-\sigma}p^{\rho}_{\lambda}(X)\circ \Phi_{\rho}$ for each vector $X$. The Levi-Civita connection $\nabla^T$ on $T(M)$ changes as 
\begin{equation}
\nabla^T{}'_XY=\nabla^T_X Y+(X\sigma)Y+(Y\sigma)X-g(X,Y)\mathrm{grad}(\sigma), 
\nonumber 
\end{equation}
where $\mathrm{grad}(\sigma):=\sum (e_i\sigma)e_i$ is the gradient vector field of $\sigma$ with respect to $g$. We show from \eqref{eqn:5-1} that the covariant derivative on $\mathbf{S}_{\rho}$ changes as 
\begin{equation}
\nabla'_X\Phi_{\rho}(\phi)-\Phi_{\rho}(\nabla_X\phi)=\Phi_{\rho}(\frac{1}{2}\sum_{ij} \{(e_i\sigma)g(X,e_j)-(e_j\sigma)g(X,e_i)\}\pi_{\rho}(e_{ij})\phi) 
\nonumber 
\end{equation}
for $\phi$ in $\Gamma(M,\mathbf{S}_{\rho})$. It follows from \eqref{eqn:4-3} and \eqref{eqn:5-14} that 
\begin{equation}
\begin{split}
&D'^{\rho}_{\lambda} \Phi_{\rho}(\phi)\\
=&\Pi^{\rho}_{\lambda}(\sum \nabla'_{e'_i}\Phi_{\rho}(\phi)\otimes (e'_i)^{\ast})=\Pi^{\rho}_{\lambda}(\sum \nabla'_{e_i}\Phi_{\rho}(\phi)\otimes (e_i)^{\ast})
=\sum p^{\rho}_{\lambda}(e_i)\nabla'_{e_i}\Phi_{\rho}(\phi)\\
=&e^{-\sigma}\Phi_{\lambda}\{D^{\rho}_{\lambda}\phi+\frac{1}{2}\sum_{ijk} p^{\rho}_{\lambda}(e_i)\{(e_k\sigma)g(e_i,e_l)-(e_l\sigma)g(e_i,e_k)\}\pi_{\rho}(e_{kl})\phi)\\
=&e^{-\sigma}\Phi_{\lambda}\{D^{\rho}_{\lambda}\phi-w(\rho;\lambda)p^{\rho}_{\lambda}(\mathrm{grad}(\sigma))\phi\}.
\end{split}\label{eqn:5-20}
\end{equation}
Here, precisely speaking, $D^{\rho}_{\lambda}$ is defined through not $\Pi^{\rho}_{\lambda}:\mathbf{S}_{\rho}\otimes T_{\mathbb{C}}(M)\to \mathbf{S}_{\lambda}$ but $\Pi^{\rho}_{\lambda}:\mathbf{S}_{\rho}\otimes T^{\ast}_{\mathbb{C}}(M)\to \mathbf{S}_{\lambda}$. So we use $\Phi_{\lambda}\circ p^{\rho}_{\lambda}(e_i)=e^{\sigma}p^{\rho}_{\lambda}(e_i)\circ \Phi_{\rho}$ in the above equation. We also have
\begin{equation}
[D^{\rho}_{\lambda}, f]=p^{\rho}_{\lambda}(\mathrm{grad}(f)) \label{eqn:5-18}
\end{equation}
for $f$ in $C^{\infty}(M)$. The equations \eqref{eqn:5-20} and \eqref{eqn:5-18} give the conformal covariance of $D^{\rho}_{\lambda}$. The next proposition answers why we call $w(\rho;\lambda)$ the conformal weight. 
\begin{prop}[\cite{F}]
When we change the Riemannian metric $g$ to $g'=e^{2\sigma}g$, the gradient $D^{\rho}_{\lambda}$ changes as 
   \begin{equation}
   D'^{\rho}_{\lambda}=(e^{(w(\rho;\lambda)-1)\sigma}\Phi_{\lambda})\circ D^{\rho}_{\lambda}\circ(e^{w(\rho;\lambda)\sigma}\Phi_{\rho})^{-1}.
    \nonumber 
   \end{equation}
   In particular, if the dimension of $\ker D^{\rho}_{\lambda}$ is finite, then $\dim \ker D^{\rho}_{\lambda}$ is a conformal invariant of $M$.  
\end{prop}

%%%%%%%%%%%%%%%%%%%%%%%%%%%%%%%%
%%           6                %%
%%%%%%%%%%%%%%%%%%%%%%%%%%%%%%%%
\section{Curvature endomorphisms}\label{sec:6}
Let $R_T$ be the Riemannian curvature $R_T$ on $T(M)$. For a local oriented orthonormal frame $e=(e_1,\cdots,e_n)$, we set a local expression of $R_T$ as $R_{ijkl}:=g(R_{T}(e_i,e_j)e_k,e_l)$, and denote the Ricci tensor by $R_{ij}=\sum_k R_{ikkj}$ and the scalar curvature by $\kappa=\sum_i R_{ii}$. We decompose the Riemannian curvature $R_{ijkl}$, 
 \begin{equation}
R_{ijkl}=W_{ijkl}+K_{ijkl}+S_{ijkl}, \nonumber 
 \end{equation}
 where 
\begin{equation}
\begin{split}
S_{ijkl}:&=\frac{\kappa}{n(n-1)}(\delta_{il}\delta_{jk}-\delta_{ik}\delta_{jl}), \\
E_{ij}:&=\frac{1}{n-2}(\frac{\kappa}{n}\delta_{ij}-R_{ij}), \\
K_{ijkl}:&=E_{ik}\delta_{jl}+E_{jl}\delta_{ik}-E_{il}\delta_{jk}-E_{jk}\delta_{il},  \\
 W_{ijkl}:&=R_{ijkl}-E_{ijkl}-S_{ijkl}.
\end{split}\label{eqn:6-6}
\end{equation}
The conformal Weyl tensor $W_{ijkl}$ and the Einstein tensor $E_{ij}$ satisfy
\begin{equation}
\sum_i W_{ijil}=0, \quad E_{ij}=E_{ji},\quad \sum_i E_{ii}=0. \label{eqn:6-7}
\end{equation}

We shall discuss curvature endomorphisms on associate vector bundle $\mathbf{S}_{\rho}$. We define the second order derivative $\nabla^2_{X,Y}$ on $\mathbf{S}_{\rho}$ for vector fields $X$ and $Y$ by 
\begin{equation}
\nabla^2_{X,Y}:=\nabla_X\nabla_Y-\nabla_{\nabla^T_X Y}. \nonumber 
\end{equation}
Then the curvature on $\mathbf{S}_{\rho}$ is $R_{\rho}(X,Y)=\nabla^2_{X,Y}-\nabla^2_{Y,X}$ for $X$ and $Y$. From \eqref{eqn:5-1}, a local expression of $R_{\rho}$ is 
\begin{equation}
R_{\rho}(e_i,e_j)=\frac{1}{2}\sum_{ij} R_{ijkl}\pi_{\rho}(e_{kl}). 
  \nonumber 
\end{equation}
By an easy calculation, we can decompose the curvature $R_{\rho}$,
\begin{equation}
R_{\rho}(e_i,e_j)=\frac{1}{2}\sum_{kl} W_{ijkl}\pi_{\rho}(e_{kl})+\sum_k(E_{ik}\pi_{\rho}(e_{kj})-E_{jk}\pi_{\rho}(e_{ki}))-\frac{\kappa}{n(n-1)}\pi_{\rho}(e_{ij}). \label{eqn:6-8}
\end{equation}

\begin{defini}
We define curvature endomorphisms in $\Gamma(M,\End(\mathbf{S}_{\rho}))$ by
\begin{equation}
R_{\rho}^q:=\sum_{ij}\pi_{\rho}(e_{ij}^q)R_{\rho}(e_i,e_j)
\nonumber 
\end{equation}
for each $q$, and when $n$ is even, 
\begin{equation}
R^{\pf}_{\rho}:=\sum_{ij} \pi_{\rho}(\pf_{ij})R_{\rho}(e_i,e_j).
\nonumber 
\end{equation}
Instead of $R_{\rho}^q$, we often use the translated curvature endomorphism
\begin{equation}
\hat{R}_{\rho}^q:=\sum_{ij}\pi_{\rho}(\hat{e}_{ij}^q)R_{\rho}(e_i,e_j)=\sum_{0\le p\le q} \binom{q}{p}\left(\frac{n-1}{2}\right)^{q-p}R^p_{\rho}. 
\nonumber 
\end{equation}
\end{defini}
\begin{ex}
We can show that $R^1_{((1/2)_{m-1},\pm 1/2)}$ is $\kappa/4$, and $R^1_{(1)}/2$ is the Ricci transformation.
\end{ex}
By using Clifford homomorphisms, the curvature endomorphisms are rewritten  as 
\begin{gather}
R_{\rho}^q=\sum_{ij,\lambda}w(\rho;\lambda)^qp^{\rho}_{\lambda}(e_i)^{\ast}p^{\rho}_{\lambda}(e_j)R_{\rho}(e_i,e_j), 
\nonumber 
 \\
R^{\pf}_{\rho}=\sum_{ij,\lambda}\pi_{\lambda}(\pf)p^{\rho}_{\lambda}(e_i)^{\ast}p^{\rho}_{\lambda}(e_j)R_{\rho}(e_i,e_j). \nonumber 
\end{gather}

\begin{prop}
The curvature endomorphisms $R_{\rho}^q$ and $R^{\pf}_{\rho}$ are self-adjoint endomorphisms of $\mathbf{S}_{\rho}$. 
\end{prop}
\begin{proof}
We consider a curvature endomorphism of $\mathbf{S}_{\rho}$, 
$$
R_{\rho}(\lambda):=\sum_{ij} p^{\rho}_{\lambda}(e_i)^{\ast}p^{\rho}_{\lambda}(e_j)R_{\rho}(e_i,e_j)=\sum_{ijkl} R_{ijkl}p^{\rho}_{\lambda}(e_i)^{\ast}p^{\rho}_{\lambda}(e_j)\pi_{\rho}(e_{kl})
$$
for each $\lambda$. From \eqref{eqn:4-32}, we have
\begin{equation}
 \begin{split}
 &\pi_{\rho}(e_{kl})p^{\rho}_{\lambda}(e_i)^{\ast}p^{\rho}_{\lambda}(e_j)-p^{\rho}_{\lambda}(e_i)^{\ast}p^{\rho}_{\lambda}(e_j)\pi_{\rho}(e_{kl})\\
=&p^{\rho}_{\lambda}(\delta_{ki}e_l-\delta_{li}e_k)p^{\rho}_{\lambda}(e_j)+p^{\rho}_{\lambda}(e_i)^{\ast}p^{\rho}_{\lambda}(\delta_{kj}e_l-\delta_{lj}e_k).
\end{split}\nonumber
\end{equation}
Then it is easy to show that 
\begin{equation}
R_{\rho}(\lambda)=\sum R_{ijkl}p^{\rho}_{\lambda}(e_i)^{\ast}p^{\rho}_{\lambda}(e_j)\pi_{\rho}(e_{kl})=\sum R_{ijkl}\pi_{\rho}(e_{kl})p^{\rho}_{\lambda}(e_i)^{\ast}p^{\rho}_{\lambda}(e_j)=R_{\rho}(\lambda)^{\ast}.\nonumber
\end{equation}
Since $R_{\rho}^q$ and $R^{\pf}_{\rho}$ are linear combinations of $\{R_{\rho}(\lambda)\}_{\lambda}$ with real coefficients, $R_{\rho}^q$ and $R^{\pf}_{\rho}$ are self-adjoint endomorphisms.
\end{proof}
We decompose the curvature endomorphisms along \eqref{eqn:6-8}. For example, we calculate the part of $\hat{R}_{\rho}^q$ depending on $E_{ij}$. From \eqref{eqn:3-16} and \eqref{eqn:6-7}, we have 
\begin{equation}
\sum_{ijk}\pi_{\rho}(\hat{e}_{ij}^q)(E_{ik}\pi_{\rho}(e_{kj})-E_{jk}\pi_{\rho}(e_{ki})) =-\sum_{ij} E_{ij}\pi_{\rho}(2\hat{e}_{ij}^{q+1}+\hat{e}_{ij}^q).
\nonumber
\end{equation}
Thus we have the decompositions of $R^q_{\rho}$ and $\hat{R}_{\rho}^q$, 
\begin{gather}
R^q_{\rho}=\frac{1}{2}\sum_{ijkl} W_{ijkl}\pi_{\rho}(e_{ij}^qe_{kl})-\sum_{ij} E_{ij}\pi_{\rho}(2e_{ij}^{q+1}+n e_{ij}^q)+\frac{\pi_{\rho}(c_{q+1})\kappa}{n(n-1)}, \nonumber\\
\hat{R}^q_{\rho}=\frac{1}{2}\sum_{ijkl} W_{ijkl}\pi_{\rho}(\hat{e}_{ij}^q\hat{e}_{kl})-\sum_{ij} E_{ij}\pi_{\rho}(2\hat{e}_{ij}^{q+1}+\hat{e}_{ij}^q)+\frac{\pi_{\rho}(2\hat{c}_{q+1}-(n-1)\hat{c}_q)\kappa}{n(n-1)}.\nonumber
\end{gather}
\begin{ex}
If $M$ is the standard sphere $S^n$, then $R^q_{\rho}$ is a constant $\pi_{\rho}(c_{q+1})$. 
\end{ex}

We consider $R^{\pf}_{\rho}$. It follows from \eqref{eqn:4-3}, \eqref{eqn:4-29}, and  \eqref{eqn:4-30} that the part of $R^{\pf}_{\rho}$ depending on the Einstein tensor is 
\begin{equation}
 \sum_{\lambda,ijk} \pi_{\lambda}(\pf)p^{\rho}_{\lambda}(e_i)^{\ast}p^{\rho}_{\lambda}(e_j)(E_{ik}\pi_{\rho}(e_{kj})-E_{jk}\pi_{\rho}(e_{ki})) =0,
\nonumber
\end{equation}
and the part depending on the scalar curvature is 
\begin{equation}
\sum_{\lambda,i,j} \pi_{\lambda}(\pf)p^{\rho}_{\lambda}(e_i)^{\ast}p^{\rho}_{\lambda}(e_j)(-\frac{\kappa}{n(n-1)}\pi_{\rho}(e_{ij}))=\frac{\pi_{\rho}(\pf)\kappa}{n-1}.\nonumber
\end{equation}
As a result, we get the following proposition which induces some interesting vanishing theorems in Section \ref{sec:8}. 
\begin{prop}\label{prop:6-3}
The curvature endomorphism $R^{\pf}_{\rho}$ does not depend on the Einstein tensor, 
\begin{equation}
R^{\pf}_{\rho}=\frac{1}{2}\sum W_{ijkl}\pi_{\rho}(\pf_{ij}e_{kl})+\frac{\pi_{\rho}(\pf)\kappa}{n-1}. \nonumber
\end{equation}
\end{prop}
\begin{cor}
\begin{enumerate}
	\item 
For the exceptional case that $\rho^{m-1}> 0$ and $\rho^{m}=0$, 
\begin{equation}
R^{\pf}_{\rho}=\frac{1}{2}\sum W_{ijkl}\pi_{\rho}(\pf_{ij}e_{kl}).\nonumber
\end{equation}
\item Suppose that $M$ is an even dimensional conformally flat manifold, then 
\begin{equation}
R^{\pf}_{\rho}=\frac{\pi_{\rho}(\pf)\kappa}{n-1}. \nonumber
\end{equation}
\end{enumerate}
\end{cor}

%%%%%%%%%%%%%%%%%%%%%%%%%%%%%%%%
%             7                 %
%%%%%%%%%%%%%%%%%%%%%%%%%%%%%%%%%
\section{Bochner-Weitzenb\"ock formulas}\label{sec:7}
The second order differential operator $(D^{\rho}_{\lambda})^{\ast}D^{\rho}_{\lambda}$ on $\mathbf{S}_{\rho}$ is realized as 
\begin{equation}
(D^{\rho}_{\lambda})^{\ast}D^{\rho}_{\lambda}=-\sum_{i,j} p^{\rho}_{\lambda}(e_i)^{\ast}p^{\rho}_{\lambda}(e_j)\nabla^2_{e_i,e_j}.  \nonumber
\end{equation}
From \eqref{eqn:4-20}, we have 
\begin{equation}
\sum_{\lambda} (D^{\rho}_{\lambda})^{\ast}D^{\rho}_{\lambda}=-\sum_{\lambda,i,j} p^{\rho}_{\lambda}(e_i)^{\ast}p^{\rho}_{\lambda}(e_j)\nabla^2_{e_i,e_j}=-\sum_{i,j} \delta_{ij}\nabla^2_{e_i,e_j}=\nabla^{\ast}\nabla, \nonumber
\end{equation}
where $\nabla^{\ast}\nabla$ is the connection Laplacian on $\mathbf{S}_{\rho}$ defined by $-\sum_i \nabla^2_{e_i,e_i}$. 

The universal Bochner-Weitzenb\"ock formula \eqref{eqn:3-18} gives
\begin{equation}
\begin{split}
\hat{R}_{\rho}^{2q}&=\sum_{i,j} \pi_{\rho}(\hat{e}_{ij}^{2q})(\nabla_{e_i,e_j}^2-\nabla_{e_j,e_i}^2) \\
 &=-\sum_{\lambda} \hat{w}(\rho;\lambda)^{2q}(D^{\rho}_{\lambda})^{\ast}D^{\rho}_{\lambda}-\sum_{i,j} \pi_{\rho}(\hat{e}_{ji}^{2q}+\sum_{p=0}^{2q-1} (-1)^p\hat{c}_{2q-1-p}\hat{e}_{ji}^p)\nabla^2_{e_j,e_i}\\
 &=\sum_{\lambda}\{\sum_{p=0}^{2q-1}\pi_{\rho}(\hat{c}_{2q-1-p})(-\hat{w}(\rho;\lambda))^p\}(D^{\rho}_{\lambda})^{\ast}D^{\rho}_{\lambda}.
\end{split}\nonumber
\end{equation}
Similaly \eqref{eqn:4-29} gives 
\begin{equation}
R^{\pf}_{\rho}=\sum_{\lambda} 2(\pi_{\rho}(\pf)-\pi_{\lambda}(\pf))(D^{\rho}_{\lambda})^{\ast}D^{\rho}_{\lambda}. 
\nonumber
\end{equation}
We are now in a position to state Bochner-Weitzenb\"ock formulas. 
\begin{thm}[Bochner-Weitzenb\"ock formulas]\label{thm:7-1}Let $\{D^{\rho}_{\lambda}\}_{\lambda}$ be the gradients on $\mathbf{S}_{\rho}$, and $\{(D^{\rho}_{\lambda})^{\ast}\}_{\lambda}$ be their formal adjoints. There exist the following identities.
\begin{gather}
\sum_{\lambda} (D^{\rho}_{\lambda})^{\ast}D^{\rho}_{\lambda}=\nabla^{\ast}\nabla,\label{eqn:7-3} \\
\sum_{\lambda} \{\sum_{p=0}^{2q-1}\pi_{\rho}(\hat{c}_{2q-1-p}) (-\hat{w}(\rho;\lambda))^p \}(D^{\rho}_{\lambda})^{\ast}D^{\rho}_{\lambda}=\hat{R}_{\rho}^{2q}, \quad q=1,2,\cdots. \label{eqn:7-4} 
\end{gather}
When $n$ is even, we also have
\begin{equation}
\sum_{\lambda} 2(\pi_{\rho}(\pf)-\pi_{\lambda}(\pf))(D^{\rho}_{\lambda})^{\ast}D^{\rho}_{\lambda} =R^{\pf}_{\rho}. \label{eqn:7-5}
\end{equation}
In the exceptional case, we set $\lambda_{\pm}:=\rho\pm \mu_m$ and have 
\begin{equation}
(D^{\rho}_{\lambda_+})^{\ast}D^{\rho}_{\lambda_+}-(D^{\rho}_{\lambda_-})^{\ast}D^{\rho}_{\lambda_-}=-\frac{1}{4\pi_{\lambda_+}(\pf)}\sum_{ijkl} W_{ijkl}\pi_{\rho}(\pf_{ij}e_{kl}). \label{eqn:7-6}
\end{equation}
\end{thm}
\begin{rem}
By using \eqref{eqn:3-19}, we obtain other identities
\begin{equation}
\hat{R}^{2q+1}_{\rho}=-\sum \{ 2\hat{w}(\rho;\lambda)^{2q+1}+\hat{w}(\rho;\lambda)^{2q}-\sum_{p=0}^{2q}\pi_{\rho}(\hat{c}_{2q-p})(-\hat{w}(\rho;\lambda))^p \} 
   (D^{\rho}_{\lambda})^{\ast}D^{\rho}_{\lambda}\nonumber 
\end{equation}
for $q=0,1,\cdots$. But, from the discussion below, these formulas are linear dependent on \eqref{eqn:7-4}. When $q$ is zero in the above equation, we have a formula by P. Gauduchon \cite{G}, 
\begin{equation}
\begin{split}
-\frac{1}{2}R^1_{\rho}&=-\frac{\pi_{\rho}(c_2)\kappa}{2n(n-1)}+\sum E_{ik}\pi_{\rho}(e_{ik}^2)-\frac{1}{4}\sum W_{ijkl}\pi_{\rho}(e_{ij}e_{kl})\\
&=\sum_{\lambda} w(\rho;\lambda)(D^{\rho}_{\lambda})^{\ast}D^{\rho}_{\lambda}.
\end{split}\label{eqn:7-8}
\end{equation}
Note that this formula can be proved from \eqref{eqn:4-21}. 
\end{rem}

We shall discuss linear independence of our Bochner-Weitzenb\"ock formulas \eqref{eqn:7-4} and \eqref{eqn:7-5}. We assume that there are $N$ irreducible components in  $\mathbf{S}_{\rho}\otimes T_{\mathbb{C}}(M)$ and $N$ gradients $\{D^{\rho}_{\lambda_i}\}_{i=1}^N$ on $\Gamma(M,\mathbf{S}_{\rho})$. In \cite{Br2}, by using spectral resolution on the standard sphere $S^n$, T. Branson showed that there is just $[N/2]$ independent identities such that $\sum_i b_{\lambda_i}(D^{\rho}_{\lambda_i})^{\ast}D^{\rho}_{\lambda_i}$ is a curvature endomorphism. So our task is to prove that the formulas \eqref{eqn:7-4} and \eqref{eqn:7-5} give $[N/2]$ independent identities. In other words, if we define the vector $v(q)$ of coefficients on \eqref{eqn:7-4} by
\begin{equation}
v(q):=(\sum_{p=0}^{2q-1}(-1)^p\pi_{\rho}(\hat{c}_{2q-1-p}) \hat{w}(\rho;\lambda_1)^p,\cdots, \sum_{p=0}^{2q-1}(-1)^p\pi_{\rho}(\hat{c}_{2q-1-p}) \hat{w}(\rho;\lambda_N)^p),  \nonumber 
\end{equation}
then we would prove that $v(1),v(2),\cdots, v([N/2])$ are linear independent in $\mathbb{R}^N$. 
We decompose $(v(1), v(2),\cdots, v(q))$ into the product of a $q\times 2q$ matrix $C(q)$ and a $2q\times N$ matrix $W(q)$ given by
\begin{gather}
C(q):=
\begin{pmatrix}
\pi_{\rho}(\hat{c}_1) & -\pi_{\rho}(\hat{c}_0) &  0 & 0  &    \cdots & 0 &  0  \\ 
\pi_{\rho}(\hat{c}_3) & -\pi_{\rho}(\hat{c}_2) & \pi_{\rho}(\hat{c}_1)  & -\pi_{\rho}(\hat{c}_0)   & \cdots & 0& 0  \\
 \cdots   &  \cdots  &  \cdots &  \cdots & \cdots    &   \cdots & \cdots \\
 \pi_{\rho}(\hat{c}_{2q-1}) & -\pi_{\rho}(\hat{c}_{2q-2}) & \cdots&\cdots&\cdots& \pi_{\rho}(\hat{c}_1) & -\pi_{\rho}(\hat{c}_0) 
\end{pmatrix}, \nonumber\\
 W(q):=\begin{pmatrix}
1           & 1            & \cdots  & 1  \\
\hat{w}(\rho;\lambda_1)& \hat{w}(\rho;\lambda_2) &\cdots & \hat{w}(\rho;\lambda_{N}) \\
\hat{w}(\rho;\lambda_1)^2& \hat{w}(\rho;\lambda_2)^2 &\cdots & \hat{w}(\rho;\lambda_{N})^2 \\
 \cdots& \cdots& \cdots &\cdots\\
 \hat{w}(\rho;\lambda_1)^{2q-1} & \hat{w}(\rho;\lambda_2)^{2q-1} &\cdots & \hat{w}(\rho;\lambda_{N})^{2q-1} 
\end{pmatrix}. \nonumber
 \end{gather}
 Since the conformal weights are differ from each other, the rank of the matrix $(v(1),v(2),\cdots, v([N/2]))=C([N/2])W([N/2])$ is $[N/2]$ except the exceptional case. For the exceptional case, the rank of $C([N/2])W([N/2])$ is $[N/2]-1$. But, there is another formula \eqref{eqn:7-6} independent of \eqref{eqn:7-4}. Thus, if there are $N$ gradients $\{D^{\rho}_{\lambda_i}\}_{i=1}^N$, then we have $[N/2]$ independent Bochner-Weitzenb\"ock formulas. 
\begin{cor}
The formulas \eqref{eqn:7-4} and \eqref{eqn:7-5} give all Bochner-Weitzenb\"ock formulas for gradients.
\end{cor}

%%%%%%%%%%%%%%%%%%%%%%%%%%%%%%%%
%%            8               %%
%%%%%%%%%%%%%%%%%%%%%%%%%%%%%%%%
\section{Examples}\label{sec:8}
We give examples and applications of our Bochner-Weitzenb\"ock formulas. We denote $D_{\lambda_i}^{\rho}$ by $D_i$ simply in this section.

We begin with fundamental examples, spinors and differential forms.
\begin{ex}[spinors] We discuss gradients on spinors, where $\Delta=((1/2)_{m})$ and  $n=2m+1$. The same result is valid for $n=2m$. We have two irreducible components in $\mathbf{S}_{\Delta}\otimes T_{\mathbb{C}}(M)$ whose highes weights are $\lambda_1=(3/2,(1/2)_{m-1})$ and $\lambda_2=\Delta=((1/2)_m)$. From \eqref{eqn:7-3} and \eqref{eqn:7-8}, we have Bochner-Weitzenb\"ock formulas
\begin{equation}
D_1^{\ast}D_1+D_2^{\ast}D_2=\nabla^{\ast}\nabla, \quad 
\frac{1}{2}D_1^{\ast}D_1-\frac{n-1}{2}D_2^{\ast}D_2=-\frac{1}{2}R_{\Delta}^1=-\frac{\kappa}{8}. \nonumber
\end{equation}
Since the Dirac operator $D$ is $\sqrt{n}D_2$ and the twistor operator $T$ is $\sqrt{n/(n-1)}D_1$, we have 
\begin{equation}
D^2=\nabla^{\ast}\nabla+\frac{1}{4}\kappa, \quad D^2=\frac{n}{4(n-1)}\kappa+T^{\ast}T. \nonumber
\end{equation}
The first equation gives a vanishing theorem. If $(M,g)$ is a compact spin manifold with positive scalar curvature, then the kernel of $D$ is zero. The second equation gives Friedrich's estimate for eigenvalues of $D^2$. Since $T^{\ast}T$ is a non-negative operator on compact spin manifold, each eigenvalue $\mu$ of $D^2$ satisfies $\mu \ge \frac{n}{4(n-1)}\min_{x\in M} \kappa(x)$ \cite{Fr}. 
\end{ex}
\begin{ex}[differential forms]
We consider the bundle of differential forms $\mathbf{S}_{(1_p)}=\Lambda^p(M)\otimes\mathbb{C} \simeq \Lambda^{n-p}(M)\otimes\mathbb{C}$ for $0\le p\le [n/2]$. We have three irreducible components in $\mathbf{S}_{(1_p)}\otimes T_{\mathbb{C}}(M)$, whose highest weights are
\[ 
\lambda_1=(2,1_{p-1}), \quad \lambda_2=(1_{p+1}),\quad \lambda_3=(1_{p-1}). 
 \]
From \eqref{eqn:7-3} and \eqref{eqn:7-8}, we obtain 
\begin{gather}
D_1^{\ast}D_1+D_2^{\ast}D_2+D_3^{\ast}D_3=\nabla^{\ast}\nabla, \nonumber \\
D_1^{\ast}D_1-p D_2^{\ast}D_2-(n-p)D_3^{\ast}D_3=-\frac{1}{2}R^1_{(1_p)}
 \nonumber
\end{gather}
The operators $D_1$, $D_2$ and $D_3$ are constant multiples of the conformal Killing operator $C$, the exterior derivative $d$ and the interior derivative $d^{\ast}$, respectively. We normalize $\{D_i\}_{1\le i\le 3}$ and obtain
\begin{gather}
C^{\ast}C+\frac{1}{p+1}d^{\ast}d+\frac{1}{n-p+1}dd^{\ast}=\nabla^{\ast}\nabla,
\nonumber 
\\
C^{\ast}C-\frac{p}{p+1}d^{\ast}d-\frac{n-p}{n-p+1}dd^{\ast}=-\frac{1}{2}R^1_{(1_p)}.\label{eqn:8-6}
\end{gather}
In particular, we get a well-known formula
\begin{equation}
d^{\ast}d+d^{\ast}d=\nabla^{\ast}\nabla+\frac{1}{2}R_{(1_p)}^1. \nonumber
\end{equation}
By using \eqref{eqn:8-6}, we can prove eigenvalue estimates of the Laplace operator $d^{\ast}d+dd^{\ast}$, Lichnerowicz's estimate for functions and Gallot-Meyer's estimate for differential forms \cite{GM}. Suppose that $(M,g)$ is a compact Riemannian manifold. For eigenfunction $f$ of $d^{\ast}d$ with nonzero eigenvalue $\mu$, we have 
\begin{equation}
\begin{split}
\mu \|df\|^2&=( dd^{\ast}df,df)=\frac{n}{n-1}( C^{\ast}C df+\frac{1}{2}R^1_{\Lambda^1}df,df)  \\
&=\frac{n}{n-1}\{\|C df\|^2+(Ric (df),df)
\}\ge \frac{n}{n-1}( Ric (df), df),
\end{split}\nonumber
\end{equation}
where $(\phi,\psi)$ denotes $\int_M\la \phi,\psi\ra dv$. Accordingly, nonzero eigenvalue of $d^{\ast}d$ on $\Lambda^0(M)$ has a lower bound depending on the Ricci curvature. 

Suppose that $(M,g)$ is a compact Riemannian manifold of positive curvature. In other words,  there exists a constant $r>0$ such that $R_{ijkl}\ge r(\delta_{il}\delta_{jk}-\delta_{ik}\delta_{jl})$. For $\phi$ in $\Gamma(M,\Lambda^p(M))$, 
\begin{equation}
\begin{split}
( (d^{\ast}d+d^{\ast}d)\phi,\phi) &\ge \frac{n-p+1}{n-p}( \{\frac{p}{p+1}d^{\ast}d+\frac{n-p}{n-p+1}dd^{\ast}\}\phi,\phi) \\
&=\frac{n-p+1}{n-p}\|C\phi\|^2+\frac{n-p+1}{2(n-p)}( R_{(1_p)}^1\phi,\phi) \\
 &\ge \frac{n-p+1}{2(n-p)}r\pi_{(1_p)}(c_2)\|\phi\|^2=p(n-p+1)r\|\phi\|^2.
\end{split}\nonumber
\end{equation}
Therefore the eigenvalue $\mu$ of $d^{\ast}d+dd^{\ast}$ on $\Lambda^p(M)$ satisfies $\mu\ge p(n-p+1)r$. 
\end{ex}
The above fundamental examples give a way of vanishing theorems and eigenvalue estimates  for gradients. We shall give more examples in this way. 
\begin{ex}[The exceptional case]
We consider the exceptional case that $n=2m$, $\rho^{m-1}>0$, and $\rho^{m}=0$. Setting $\lambda_{\pm}:=\rho\pm \mu_m$ and $D_{\pm}:=D^{\rho}_{\lambda_{\pm}}$, we have 
\begin{equation}
D_+^{\ast}D_+-D_-^{\ast}D_-=-\frac{1}{4\pi_{\lambda_+}(\pf)}\sum_{ijkl} W_{ijkl}\pi_{\rho}(\pf_{ij}e_{kl}). \nonumber
\end{equation}
If $(M,g)$ is a conformally flat manifold, then $D_+^{\ast}D_+=D_-^{\ast}D_-$. When $\rho$ is $(1_{m-1})$, $D_{\pm}$ is a constant multiple of $d^{\pm}=1/2(1\pm \ast)d$ on $\Gamma(M,\Lambda^{m-1}(M))$ and $D_+^{\ast}D_+=D_-^{\ast}D_-$ automatically. 
\end{ex}
\begin{ex}[A vanishing theorem associated to the Pfaffian element]
We consider the case of $n=2m$ and $\rho=(p_{m-1},\pm p)$ for $p=1,3/2,2,5/2\cdots$. We have two irreducible components in $\mathbf{S}_{\rho}\otimes T_{\mathbb{C}}(M)$ and set $\lambda_1=(p+1,p_{m-2},\pm p)$ and $\lambda_2=(p_{m-1},\pm p \mp 1)$. Then we have 
\begin{gather}
D_1^{\ast}D_1+D_2^{\ast}D_2=\nabla^{\ast}\nabla, \nonumber \\
p D_1^{\ast}D_1+(-p-m+1)D_2^{\ast}D_2=-\frac{1}{2}R^1_{\rho},\nonumber\\
2(\pi_{\rho}(\pf)-\pi_{\lambda_1}(\pf))D_1^{\ast}D_1+2(\pi_{\rho}(\pf)-\pi_{\lambda_2}(\pf))D_2^{\ast}D_2=R^{\pf}_{\rho}.\nonumber
\end{gather}
The second identity is linear dependent on the third one. In fact, we have $\pi_{\rho}(\pf)R^1_{\rho}=p(p+m-1)R^{\pf}_{\rho}$. Then $R^1_{\rho}$ does not depend on the Einstein tensor,
\begin{equation}
R^1_{\rho}=\frac{p(p+m-1)}{2m-1}\kappa+\frac{1}{2}\sum W_{ijkl}\pi_{\rho}(e_{ij}e_{kl}). \nonumber
\end{equation}
Thus we have
\begin{equation}
\frac{2p+m-1}{p}D_2^{\ast}D_2=\nabla^{\ast}\nabla+\frac{p+m-1}{4m-2}\kappa+\frac{1}{4p}\sum W_{ijkl}\pi_{\rho}(e_{ij}e_{kl}). \nonumber
\end{equation}
We conclude that, if $M$ is a compact conformally flat manifold with positive scalar curvature, then the kernel of $D_2$ is zero. In particular, for $\rho=(1_{m-1},\pm 1)$, we have Bourguignon's vanishing $H^{m}(M,\mathbb{R})=0$ \cite{Bo}. This vanishing theorem for $D_2$ has been shown by T. Branson and O. Hijazi \cite{BH1}, \cite{BH2}. They also discussed a relation to the first eigenvalue of Yamabe Laplacian. 

We shall generalize the above example. We consider the case that $n=2m$ and $\rho^m\neq 0$, and order highest weights as $\lambda_1=\rho+\mu_1>\lambda_2>\cdots >\lambda_N$.  It follows from \eqref{eqn:7-3} and \eqref{eqn:7-5} that 
\begin{equation}
\nabla^{\ast}\nabla+\frac{1}{2(\pi_{\lambda_1}(\pf)-\pi_{\rho}(\pf))}R^{\pf}_{\rho}=\sum_{2\le i \le N}\frac{\pi_{\lambda_1}(\pf)-\pi_{\lambda_i}(\pf)}{\pi_{\lambda_1}(\pf)-\pi_{\rho}(\pf)}D_i^{\ast}D_i, \nonumber
\end{equation}
where the coefficient of $D_{i}^{\ast}D_i$ is positive. If $M$ is a compact conformally flat manifold with positive scalar curvature, then $\bigcap_{2\le i\le N}\ker D_i=0$. 
\end{ex}
\begin{ex}[The conformal Weyl tensor]
Let $(M,g)$ be an $n$-dimensional Riemannian manifold with $n\ge 5$. The four dimensional case is discussed in the next section. The conformal Weyl tensor $W=W_{ijkl}$ is a section of $\mathbf{S}_{\rho}$ with $\rho=(2_2)$. The highest weights of irreducible components of $\mathbf{S}_{\rho}\otimes T_{\mathbb{C}}(M)$ are \begin{equation}
\begin{cases}
(3,2), \quad (2_2) \quad (2,1),  & n=5, \\
(3,2), \quad (2_2,1),\quad (2_2,-1),\quad (2,1), &n=6, \\
(3,2) \quad (2_2,1),\quad (2,1), & n\ge 7.
\end{cases}\nonumber
\end{equation}
Because of the second Bianchi identities, the projections of $\nabla W$ to $\mathbf{S}_{(2_2,\pm 1)}$ and $\mathbf{S}_{(2_2)}$ are zero \cite{S}. So we set $\lambda_1:=(3,2)$ and $\lambda_2:=(2,1)$. It follows from \eqref{eqn:7-3} and \eqref{eqn:7-8} that there is a Bochner-Weitzenb\"ock formula for the conformal Weyl tensor,
\begin{equation}
\frac{n+1}{2}D_{2}^{\ast}D_{2}W=\nabla^{\ast}\nabla W+\frac{1}{4}R^1_{\rho}W.
\nonumber
\end{equation}
It is known that $D_2 W$ is a constant multiple of $\sum_s (\nabla^s W_{sijk}+\nabla^s W_{sjik})$ \cite{S}. Then $D_{2}W=0$ is equivalent to $\delta W=0$, where $\delta W:=-\sum \nabla^s W_{sijk}$. If $(M,g)$ is a Riemannian manifold with $\delta W=0$, then $W$ satisfies
\begin{equation}
\nabla^{\ast}\nabla W+\frac{1}{4}R^1_{\rho}W=0.\nonumber
\end{equation}
\end{ex}
%%%%%%%%%%%%%%%%%%%%%%%%%%%%%%%%%%%
%%%%%%%%%%%%%%%%%%%%%%%%%%%%%%%%
%%            9               %%
%%%%%%%%%%%%%%%%%%%%%%%%%%%%%%%%
\section{The four dimensional case}\label{sec:9}
In this section, we discuss the four dimensional case. We begin with basic facts of four dimensional geometry. The Hodge star operator $\ast$ decomposes $\Lambda^2(\mathbb{R}^4)$ into self-dual part $\Lambda^2_+$ and anti-self-dual part $\Lambda^2_-$. For an oriented orthonormal basis $e=(e_1,e_2,e_3,e_4)$ of $\mathbb{R}^4$, we set a basis of $\Lambda^2_{\pm}$ as
\begin{equation}
X_1^{\pm}:=\frac{1}{2}(e_{14}\pm e_{23}), \quad X_2^{\pm}:=-\frac{1}{2}(e_{13}\pm e_{42}),\quad X_3^{\pm}:=\frac{1}{2}(e_{12}\pm e_{34}).\nonumber
\end{equation}
Identifying  $\Lambda^2(\mathbb{R}^4)$ with $\mathfrak{so}(4)=\mathfrak{so}(3)\oplus \mathfrak{so}(3)$, we know 
\begin{equation}
[X_i^{\pm},X_j^{\pm}]=\sum_{1\le k\le 3} \epsilon_{ijk}X_k^{\pm},\quad [X_i^+,X_j^-]=0,
\nonumber
\end{equation}
where 
\begin{equation}
\epsilon_{ijk}=\begin{cases}
\mathrm{sgn}\left(
\begin{smallmatrix}
 1& 2 & 3\\ 
i & j& k
\end{smallmatrix}\right),  &  \{i,j,k\}=\{1,2,3\}, \\
0, & \textrm{otherwise}. 
\end{cases}\nonumber
\end{equation}
Thus $\{X_i^{\pm}\}_i$ is a standard basis of $\mathfrak{so}(3)$, 
\begin{equation}
X_1^{\pm}=\begin{pmatrix}
0 & 0 & 0\\
0 & 0& -1 \\
0 & 1& 0
\end{pmatrix},\quad 
X_2^{\pm}=\begin{pmatrix}
0 & 0 & 1\\
0 & 0& 0 \\
-1 & 0& 0
\end{pmatrix},\quad 
X_3^{\pm}=\begin{pmatrix}
0 & -1 & 0\\
1 & 0& 0 \\
0 & 0& 0
\end{pmatrix}. \label{eqn:9-4}
\end{equation}

Let $(M,g)$ be a four dimensional oriented Riemannian or spin manifold. We think of the Riemannian curvature $R_T$ as an endomorphism of $\Lambda^2(M)$, 
\begin{equation}
R_T:\Lambda^2(M)\ni e_{ij}\mapsto \frac{1}{2}\sum R_{ijkl}e_{kl}\in \Lambda^2(M). \nonumber
\end{equation}
Then we realize $R_T$ as a $6\times 6$ matrix with respect to basis $\{X_i^+\}_i\cup \{X_j^-\}_j$, 
\begin{equation}
\begin{pmatrix}
W^+ & 0 \\
0 & W^-
\end{pmatrix}+
\begin{pmatrix}
0 & K \\
K^t & 0
\end{pmatrix}+
\begin{pmatrix}
-\kappa/12 &0 \\
0 & -\kappa/12
\end{pmatrix},\nonumber
\end{equation}
where $W^+$ (resp. $W^-$) is the self-dual (resp. anti-self-dual) conformal Weyl tensor and $K$ corresponds to $K_{ijkl}$ in \eqref{eqn:6-6}. 
In other words, we have 
\begin{equation}
\begin{split}
R_T(X_i^{+})=\sum W^{+}_{ij}X_j^{+}+\sum K_{ji}X_j^{-}-\frac{\kappa}{12}X_i^{+}
, \\
R_T(X_i^-)=\sum W^-_{ij}X^-_j+\sum K_{ij}X_j^+-\frac{\kappa}{12}X_i^-.
\end{split}\nonumber
\end{equation}

Now we denote the highest weight of $\mathfrak{so}(4)$ by $\rho=(\frac{k+l}{2},\frac{k-l}{2})$ for non-negative integers $k$ and $l$. Note that $\rho$ corresponds to  $(k/2)\hat{\otimes}(l/2)$ as a highest weight of $\mathfrak{so}(3)\oplus \mathfrak{so}(3)$.

We define curvature endomorphisms $R^{\pm}_{\rho}$ on $\mathbf{S}_{\rho}$ by 
$$
R^{\pm}_{\rho}:=4\sum_i\pi_{\rho}(X_i^{\pm}R_T(X_i^{\pm})).
$$
Since $4\sum X_i^+X_i^+=-c_2/2-2\pf$ and $4\sum X_i^- X_i^-=-c_2/2+2\pf$, we know 
\begin{equation}
\begin{split}
R^{+}_{\rho}&=4\sum_{ij} W_{ij}^+\pi_{\rho}(X_i^{+}X_j^{+})+4\sum_{ij} K_{ij}\pi_{\rho}(X_i^{-}X_j^{+})+\frac{k(k+2)}{12}\kappa,\\
   R^{-}_{\rho}&=4\sum_{ij} W_{ij}^{-}\pi_{\rho}(X_i^{-}X_j^{-})+4\sum_{ij} K_{ij}\pi_{\rho}(X_i^{-}X_j^{+})+\frac{l(l+2)}{12}\kappa.
   \end{split}\nonumber
\end{equation}

We state some properties of these curvature endomorphisms.
\begin{lem}
\begin{enumerate}
\item The curvature endomorphisms $R_{\rho}^1$ and $R_{\rho}^{\pf}$ are realized as 
\begin{equation}
R_{\rho}^1=R_{\rho}^++R_{\rho}^-, \quad R_{\rho}^{\pf}=R_{\rho}^+-R_{\rho}^-.
 \label{eqn:9-9}
\end{equation}
\item If  $\rho$ is $(k/2,k/2)$, then  $R^{-}_{\rho}$ is zero and $R^+_{\rho}$ does not depend on the Einstein tensor. If $\rho$ is $(\frac{k+1}{2},\frac{k-1}{2})$, then $R^{-}_{\rho}$ does not depend the anti-self-dual conformal Weyl tensor $W^-$. 
\end{enumerate}
\end{lem}
\begin{proof}
By using \eqref{eqn:pf-4}, we can show \eqref{eqn:9-9} straightforwardly. We shall prove the second claim. When $\rho$ is $(k/2,k/2)$, $\pi_{\rho}(X_i^-)$ is zero for each $i$. Then $R^{-}_{\rho}$ is zero and $R^+_{\rho}$ does not depend on the Einstein tensor. When $\rho$ is $(\frac{k+1}{2},\frac{k-1}{2})$, $\{\pi_{\rho}(X_i^-)\}_{1\le i\le 3}$ give a spin $1/2$ representation of $\mathfrak{so}(3)=\mathfrak{su}(2)$. Then we have $\sum_{ij} W_{ij}^{-}\pi_{\rho}(X_i^{-}X_j^{-})=0$ because of $W_{ij}^-=W_{ji}^-$ and $\sum W_{ii}^-=0$. 
\end{proof}
 The irreducible decomposition of $\mathbf{S}_{\rho}\otimes T_{\mathbb{C}}(M)$ is $\oplus_{i=1}^4 \mathbf{S}_{\lambda_i}$, where 
 \begin{equation}
\lambda_1=\rho+\mu_1,\quad \lambda_2=\rho+\mu_2,\quad  \lambda_3=\rho-\mu_2, \quad  \lambda_4=\rho-\mu_1. 
\nonumber 
\end{equation}
From \eqref{eqn:7-3}, \eqref{eqn:7-4} and \eqref{eqn:7-5},  we have all Bochner-Weitzenb\"ock formulas for the four dimensional case, 
\begin{equation}
\begin{split}
D_1^{\ast}D_1+D_2^{\ast}D_2+D_3^{\ast}D_3+D_4^{\ast}D_4=\nabla^{\ast}\nabla, \\
k D_1^{\ast}D_1+k D_2^{\ast}D_2-(k+2)D_3^{\ast}D_3-(k+2)D_4^{\ast}D_4=-R_{\rho}^+, \\
l D_1^{\ast}D_1-(l+2)D_2^{\ast}D_2+l D_3^{\ast}D_3-(l+2)D_4^{\ast}D_4=-R_{\rho}^-,
\end{split}\nonumber
\end{equation}
where $D_i$ denotes $D^{\rho}_{\lambda_i}$. 

We shall state some vanishing theorems. First, we consider the case that $\rho$ is $(k/2,k/2)$ for a positive integer $k$. Then $D_2=D_4=0$ and 
\begin{equation}
\frac{2(k+1)}{k}D_3^{\ast}D_3=\nabla^{\ast}\nabla+\frac{4}{k}\sum W_{ij}^+\pi_{\rho}(X_i^+X_j^+)+\frac{k+2}{12}\kappa. \nonumber
\end{equation}
\begin{ex}
If $M$ is a compact anti-self-dual manifold, then 
\begin{equation}
\frac{2(k+1)}{k}D_3^{\ast}D_3=\nabla^{\ast}\nabla+\frac{k+2}{12}\kappa. 
    \nonumber
\end{equation}
The kernel of $D_3$ is isomorphic to a cohomology on the twistor space of $M$ \cite{Hi}. 
\end{ex}
\begin{ex}[self-dual $2$-forms]
For $\rho=(1,1)$,  the associated vector bundle $\mathbf{S}_{\rho}$ is $\Lambda^2_+(M)\otimes\mathbb{C} $ and the kernel of $D_3$ is the space of harmonic self-dual $2$-forms. If $\omega$ is a harmonic self-dual  $2$-form, then we have a formula in \cite{Bo},
\begin{equation}
0=\la\nabla^{\ast}\nabla\omega,\omega\ra+2\la W^+(\omega),\omega\ra+\frac{1}{3}\kappa\la\omega,\omega\ra, \nonumber
\end{equation}
where $W^+(\omega)=W^+(\sum \omega_{i}X_i^+)=\sum \omega_i W^+_{ij}X_j^+$. 
\end{ex}
\begin{ex}[the self-dual conformal Weyl tensor]
For $\rho=(2,2)$, the vector bundle $\mathbf{S}_{\rho}$ is embedded into $\Lambda^2_+(M)\otimes (\Lambda^2_+(M))^{\ast}\otimes\mathbb{C}$ and the  sections of $\mathbf{S}_{\rho}$ are realized locally as trace-free symmetric $3\times 3$ matrices. Therefore we have $\la Z,Z\ra=\mathrm{tr}(Z\bar{Z})$ and $\pi_{\rho}(X_i^+)Z=X_i^+Z-Z X_i^+$ for $Z$ in $\Gamma(M,\mathbf{S}_{\rho})$. Here we use the matrix realization \eqref{eqn:9-4} of $\{X_i^+\}_i$. It follows that 
\begin{equation}
\begin{split}
\sum W_{ij}^+\pi_{\rho}(X_i^+X_j^+)Z &=\sum W_{ij}^+(X_j^+X_i^+Z-X_j^+Z X_i^+-X_i^+Z X_j^++Z X_i^+X_j^+) \\
&=W^+Z+ZW^+-2\sum W_{ij}^+X_i^+Z X_j^+.
\end{split}\nonumber
\end{equation}
Choosing $W^+$  as a section $Z$ of $\mathbf{S}_{\rho}$, we have 
\begin{equation}
\begin{split}
\la\sum W_{ij}^+\pi_{\rho}(X_i^+X_j)W^+,W^+\ra&=2\mathrm{tr}((W^+)^3)+4\mathrm{tr}(W^+)^3)=6\mathrm{tr}((W^+)^3)\\
&=18\det(W^+)
\end{split} \nonumber
\end{equation}
If $\delta W^+=-\sum_s \nabla^s W^+_{sijk}=0$, then $D_3 W^+$ is zero and $W^+$ satisfies the following formula in \cite{Be},
\begin{equation}
0=\la \nabla^{\ast}\nabla W^+,W^+\ra+6\mathrm{tr}((W^+)^3)+\frac{\kappa}{2}\mathrm{tr}((W^+)^2). \nonumber
\end{equation}
\end{ex}

Next we consider the case that the highest weight $\rho$ is $(\frac{k+1}{2},\frac{k-1}{2})$ for a positive integer $k$. Then we have 
\begin{equation}
\begin{split}
\nabla^{\ast}\nabla-1/3R^-_{\rho}=4/3(D_1^{\ast}D_1+D_3^{\ast}D_3), \\
\nabla^{\ast}\nabla+R^-_{\rho}=4(D_2^{\ast}D_2+D_4^{\ast}D_4), \nonumber
\end{split}
\end{equation}
where
\begin{equation}
R^-_{\rho}=4\sum K_{ij}\pi_{\rho}(X_i^-X_j^+)+\kappa/4. \nonumber
\end{equation}
\begin{ex}
If $(M,g)$ is a compact Einstein manifold with positive (reps. negative) constant scalar curvature, then $\ker D_2\cap \ker D_4$ (resp. $\ker D_1\cap \ker D_3$) is zero. 
\end{ex}

%%%%%%%%%%%%%%%%%%%%%%%%%%%%%%%
%%% Acknowledgement %%%%%%%%%%%%
%%%%%%%%%%%%%%%%%%%%%%%%%%%%%%%
\section*{Acknowledgement}
The author is partially supported by the Grant-in-Aid for JSPS Fellows from the Ministry of Education, Culture, Sports, Science and Technology. He thanks David M. J. Calderbank for some useful comments. 

%%%%%%%%%%%%%%%%%%%%%%%%%%%%%%%%
%        reference             %
%%%%%%%%%%%%%%%%%%%%%%%%%%%%%%%%
%%

%%%%%%%%%%%%%%%%%%%%%%%%%
\begin{flushleft}
Yasushi Homma \\
Department of Mathematics, \\
 Faculty of Science and Technology, \\
 Science University of Tokyo, \\ 
 2641 Noda, Chiba, 278-8510, \\
 JAPAN. \\
  \textit{E-mail address}: homma\_yasushi@ma.noda.tus.ac.jp
\end{flushleft}

\end{document}